# Model selection for Poisson processes

**Lucien Birgé**[1]

*Université Paris VI*

**Abstract:** Our purpose in this paper is to apply the general methodology for model selection based on T-estimators developed in Birgé [*Ann. Inst. H. Poincaré Probab. Statist.* **42** (2006) 273–325] to the particular situation of the estimation of the unknown mean measure of a Poisson process. We introduce a Hellinger type distance between finite positive measures to serve as our loss function and we build suitable tests between balls (with respect to this distance) in the set of mean measures. As a consequence of the existence of such tests, given a suitable family of approximating models, we can build T-estimators for the mean measure based on this family of models and analyze their performances. We provide a number of applications to adaptive intensity estimation when the square root of the intensity belongs to various smoothness classes. We also give a method for aggregation of preliminary estimators.

## 1. Introduction

This paper deals with the estimation of the mean measure $\mu$ of a Poisson process $\boldsymbol{X}$ on $\mathcal{X}$. More precisely, we develop a theoretical, but quite general method for estimating $\mu$ by model selection with applications to adaptive estimation and aggregation of preliminary estimators. The main advantage of the method is its generality. We do not make any assumption on $\mu$ apart from the fact that it should be finite and we allow arbitrary countable families of models provided that each model be of finite metric dimension, i.e. is not too large in a suitable sense to be explained below. We do not know of any other estimation method allowing to deal with model selection in such a generality and with as few assumptions. The main drawback of the method is its theoretical nature, effective computation of the estimators being typically computationally too costly for permitting a practical implementation. In order to give a more precise idea of what this paper is about, we need to start by recalling a few well-known facts about Poisson processes that can, for instance, be found in Reiss [29].

### 1.1. The basics of Poisson processes

Let us denote by $\mathcal{Q}_+(\mathcal{X})$ the cone of finite positive measures on the measurable space $(\mathcal{X}, \mathcal{E})$. Given an element $\mu \in \mathcal{Q}_+(\mathcal{X})$, a Poisson process on $\mathcal{X}$ with mean measure $\mu$ is a point process $\boldsymbol{X} = \{X_1, \ldots, X_N\}$ on $\mathcal{X}$ such that $N$ has a Poisson distribution with parameter $\mu(\mathcal{X})$ and, conditionally on $N$, the $X_i$ are i.i.d. with distribution $\mu_1 = \mu/\mu(\mathcal{X})$. Equivalently, the Poisson process can be viewed as a random measure $\Lambda_{\boldsymbol{X}} = \sum_{i=1}^{N} \delta_{X_i}$, $\delta_x$ denoting the Dirac measure concentrated

---

at the point $x$. Then, whatever the partition $A_1, \ldots, A_n$ of $\mathcal{X}$, the $n$ random variables $\Lambda_{\boldsymbol{X}}(A_i)$ are independent with Poisson distributions and respective parameters $\mu(A_i)$ and this property characterizes a Poisson process. We shall denote by $Q_\mu$ the distribution of a Poisson process with mean measure $\mu$ on $\mathcal{X}$. We recall that, for any nonnegative measurable function $\phi$ on $(\mathcal{X}, \mathcal{E})$,

$$\text{(1.1)} \qquad \mathbb{E}\left[\sum_{i=1}^{N} \phi(X_i)\right] = \int_{\mathcal{X}} \phi(x)\, d\mu(x)$$

and

$$\text{(1.2)} \qquad \mathbb{E}\left[\prod_{i=1}^{N} \phi(X_i)\right] = \exp\left[\int_{\mathcal{X}} [\phi(x) - 1]\, d\mu(x)\right].$$

If $\mu, \nu \in \mathcal{Q}_+(\mathcal{X})$ and $\mu \ll \nu$, then $Q_\mu \ll Q_\nu$ and

$$\text{(1.3)} \qquad \frac{dQ_\mu}{dQ_\nu}(X_1, \ldots, X_N) = \exp[\nu(\mathcal{X}) - \mu(\mathcal{X})] \prod_{i=1}^{N} \frac{d\mu}{d\nu}(X_i),$$

with the convention that $\prod_{i=1}^{0} (d\mu/d\nu)(X_i) = 1$.

### 1.2. Introducing our loss function

From now on, we assume that we observe a Poisson process $\boldsymbol{X}$ on $\mathcal{X}$ with unknown mean measure $\mu \in \mathcal{Q}_+(\mathcal{X})$ so that $\mu$ always denotes the parameter to be estimated. For this, we use estimators $\hat{\mu}(\boldsymbol{X})$ with values in $\mathcal{Q}_+(\mathcal{X})$ and measure their performance via the loss function $H^q(\hat{\mu}(\boldsymbol{X}), \mu)$ for $q \geq 1$, where $H$ is a suitable distance on $\mathcal{Q}_+(\mathcal{X})$. To motivate its introduction, let us recall some known facts. The Hellinger distance $h$ between two probabilities $P$ and $Q$ defined on the same space and their Hellinger affinity $\rho$ are given respectively by

$$\text{(1.4)} \quad h^2(P, Q) = \frac{1}{2}\int \left(\sqrt{dP} - \sqrt{dQ}\right)^2, \quad \rho(P, Q) = \int \sqrt{dP\, dQ} = 1 - h^2(P, Q),$$

where $dP$ and $dQ$ denote the densities of $P$ and $Q$ with respect to any dominating measure, the result being independent of the choice of such a measure. If $X_1, \ldots, X_n$ are i.i.d. with distribution $\overline{P}$ on $\mathcal{X}$ and $\overline{Q}$ is another distribution, it follows from an exponential inequality that, for all $x \in \mathbb{R}$,

$$\text{(1.5)} \quad \begin{aligned} \mathbb{P}\left[\sum_{i=1}^{n} \log\left(\frac{d\overline{Q}}{d\overline{P}}\right)(X_i) \geq 2x\right] &\leq \exp\left[n \log\left(\rho(\overline{P}, \overline{Q})\right) - x\right] \\ &\leq \exp\left[nh^2\left(\overline{P}, \overline{Q}\right) - x\right], \end{aligned}$$

which provides an upper bound for the errors of likelihood ratio tests. In particular, if $\mu$ and $\mu'$ are two elements in $\mathcal{Q}_+(\mathcal{X})$ dominated by some measure $\lambda$, it follows from (1.3) and (1.2) that the Hellinger affinity $\rho(Q_\mu, Q_{\mu'})$ between $\mu$ and $\mu'$ is given by

$$\text{(1.6)} \qquad \rho(Q_\mu, Q_{\mu'}) = \int \sqrt{\frac{dQ_\mu}{dQ_\lambda}\frac{dQ_{\mu'}}{dQ_\lambda}}\, dQ_\lambda = \exp\left[-H^2(\mu, \mu')\right],$$

where

$$H^2(\mu, \mu') = \frac{1}{2}\left[\mu(\mathcal{X}) + \mu'(\mathcal{X})\right] - \int \sqrt{(d\mu/d\lambda)(d\mu'/d\lambda)} \tag{1.7}$$

$$= \frac{1}{2}\int \left(\sqrt{d\mu/d\lambda} - \sqrt{d\mu'/d\lambda}\right)^2. \tag{1.8}$$

Comparing (1.8) with (1.4) indicates that $H$ is merely the generalization of the Hellinger distance $h$ between probabilities to arbitrary finite positive measures and the introduction of $H$ turns $\mathcal{Q}_+(\mathcal{X})$ into a metric space. Moreover, we derive from (1.5) with $n = 1$ that, when $\boldsymbol{X}$ is a Poisson process with mean measure $\mu$ on $\mathcal{X}$,

$$\mathbb{P}\left[\log\left(\frac{dQ_{\mu'}}{dQ_\mu}\right)(\boldsymbol{X}) \geq 2x\right] \leq \exp\left[-H^2(\mu, \mu') - x\right]. \tag{1.9}$$

If $\mu(\mathcal{X}) = \mu'(\mathcal{X}) = n$, then $H^2(\mu, \mu') = nh^2(\mu_1, \mu'_1)$ and (1.9) becomes a perfect analogue of (1.5). The fact that the errors of likelihood ratio tests between two probabilities are controlled by their Hellinger affinity justifies the introduction of the Hellinger distance as the natural loss function for density estimation, as shown by Le Cam [26]. It also motivates the choice of $H^q$ as a natural loss function for estimating the mean measure of a Poisson process. For simplicity, we shall first focus on the quadratic risk $\mathbb{E}[H^2(\hat{\mu}(\boldsymbol{X}), \mu)]$.

### 1.3. Intensity estimation

A case of particular interest occurs when we have at hand a reference positive measure $\lambda$ on $\mathcal{X}$ and we assume that $\mu \ll \lambda$ with $d\mu/d\lambda = s$, in which case $s$ is called the *intensity* (with respect to $\lambda$) of the process with mean measure $\mu$. Denoting by $\mathbb{L}_i^+(\lambda)$ the positive part of $\mathbb{L}_i(\lambda)$ for $i = 1, 2$, we observe that $s \in \mathbb{L}_1^+(\lambda)$, $\sqrt{s} \in \mathbb{L}_2^+(\lambda)$ and $\mu \in \mathcal{Q}_\lambda = \{\mu_t = t \cdot \lambda, t \in \mathbb{L}_1^+(\lambda)\}$. The one-to-one correspondence $t \mapsto \mu_t$ between $\mathbb{L}_1^+(\lambda)$ and $\mathcal{Q}_\lambda$ allows us to transfer the distance $H$ to $\mathbb{L}_1^+(\lambda)$ which gives, by (1.8),

$$H(t, u) = H(\mu_t, \mu_u) = \left(1/\sqrt{2}\right) \left\|\sqrt{t} - \sqrt{u}\right\|_2 \quad \text{for } t, u \in \mathbb{L}_1^+(\lambda), \tag{1.10}$$

where $\|\cdot\|_2$ stands for the norm in $\mathbb{L}_2(\lambda)$. When $\mu = \mu_s \in \mathcal{Q}_\lambda$ it is natural to estimate it by some element $\hat{\mu}(\boldsymbol{X}) = \hat{s}(\boldsymbol{X}) \cdot \lambda$ of $\mathcal{Q}_\lambda$, in which case $H(\hat{\mu}(\boldsymbol{X}), \mu) = H(\hat{s}(\boldsymbol{X}), s)$ and our problem can be viewed as a problem of *intensity estimation*: design an estimator $\hat{s}(\boldsymbol{X}) \in \mathbb{L}_1^+(\lambda)$ for the unknown intensity $s$. From now on, given a Poisson process $\boldsymbol{X}$ with mean measure $\mu$, we shall denote by $\mathbb{E}_\mu$ and $\mathbb{P}_\mu$ (or $\mathbb{E}_s$ and $\mathbb{P}_s$ when $\mu = \mu_s$) the expectations of functions of $\boldsymbol{X}$ and probabilities of events depending on $\boldsymbol{X}$, respectively.

### 1.4. Model based estimation and model selection

It is common practice to try to estimate the intensity $s$ on $\mathcal{X}$ by a piecewise constant function, i.e. a histogram estimator $\hat{s}(\boldsymbol{X})$ belonging to the set

$$\overline{S}_m = \left\{\sum_{j=1}^D a_j \mathbb{1}_{I_j}, \, a_j \geq 0 \text{ for } 1 \leq j \leq D\right\}$$

of nonnegative piecewise constant functions with respect to the partition $\{I_1, \ldots, I_D\} = m$ of $\mathcal{X}$ with $\lambda(I_j) > 0$ for all $j$. More generally, given a finite family $m = \{\varphi_1, \ldots, \varphi_D\}$ of elements of $\mathbb{L}_2(\lambda)$, we may consider the $D$-dimensional linear space $\overline{S}_m$ generated by the $\varphi_j$ and try to estimate $\sqrt{s}$ by some element $\sqrt{\hat{s}(\boldsymbol{X})} \in \overline{S}_m$. This clearly leads to difficulties since $\overline{S}_m$ is not a subset of $\mathbb{L}_2^+(\lambda)$, but we shall nevertheless show that it is possible to design an estimator $\hat{s}_m(\boldsymbol{X})$ with the property that

$$(1.11) \qquad \mathbb{E}_s\left[H^2\left(\hat{s}_m(\boldsymbol{X}), s\right)\right] \leq C\left[\inf_{t \in \overline{S}_m} \left\|t - \sqrt{s}\right\|_2^2 + |m|\right],$$

where $|m| = D$ stands for the cardinality of $m$ and $C$ is a universal constant. In this approach, $\overline{S}_m$ should be viewed as a *model* for $\sqrt{s}$, which means an approximating set since we never assume that $\sqrt{s} \in \overline{S}_m$ and the risk bound (1.11) has (up to the constant $C$) the classical structure of the sum of an approximation term $\inf_{t \in \overline{S}_m} \|t - \sqrt{s}\|_2^2$ and an estimation term $|m|$ corresponding to the number of parameters to be estimated.

If we introduce a countable (here countable always means finite or countable) family of models $\{\overline{S}_m, m \in \mathcal{M}\}$ of the previous form, we would like to know to what extent it is possible to build a new estimator $\hat{s}(\boldsymbol{X})$ such that

$$(1.12) \qquad \mathbb{E}_s\left[H^2\left(\hat{s}(\boldsymbol{X}), s\right)\right] \leq C' \inf_{m \in \mathcal{M}} \left\{\inf_{t \in \overline{S}_m} \left\|t - \sqrt{s}\right\|_2^2 + |m|\right\},$$

for some other constant $C'$, i.e. to know whether one can design an estimator which realizes, up to some constant, the best compromise between the two components of the risk bound (1.11). The problem of understanding to what extent (1.12) does hold has been treated in many papers using various methods, mostly based on the minimization of some penalized criterion. A special construction based on testing has been introduced in Birgé [9] and then applied to different stochastic frameworks. We shall show here that this construction also applies to Poisson processes and then derive the numerous consequences of this property. We shall, in particular, be able to prove the following result in Section 3.4.1 below.

**Theorem 1.** *Let $\lambda$ be some positive measure on $\mathcal{X}$ and $\|\cdot\|_2$ denote the norm in $\mathbb{L}_2(\lambda)$. Let $\{\overline{S}_m\}_{m \in \mathcal{M}}$ be a finite or countable family of linear subspaces of $\mathbb{L}_2(\lambda)$ with respective finite dimensions $\overline{D}_m$ and let $\{\Delta_m\}_{m \in \mathcal{M}}$ be a family of nonnegative weights satisfying*

$$(1.13) \qquad \sum_{m \in \mathcal{M}} \exp[-\Delta_m] \leq \Sigma < +\infty.$$

*Let $\boldsymbol{X}$ be a Poisson process on $\mathcal{X}$ with unknown mean measure $\mu = \mu_s + \mu^\perp$ where $s \in \mathbb{L}_1^+(\lambda)$ and $\mu^\perp$ is orthogonal to $\lambda$. One can build an estimator $\hat{\mu} = \hat{\mu}(\boldsymbol{X}) = \hat{s}(\boldsymbol{X}) \cdot \lambda \in \mathcal{Q}_\lambda$ satisfying, for all $\mu \in \mathcal{Q}_+(\mathcal{X})$ and $q \geq 1$,*

$$(1.14) \quad \begin{aligned} \mathbb{E}_\mu\left[H^q(\mu, \hat{\mu})\right] &\leq C(q)\left[1 + \Sigma\right] \\ &\times \left[\sqrt{\mu^\perp(\mathcal{X})} + \inf_{m \in \mathcal{M}}\left\{\inf_{t \in \overline{S}_m} \left\|\sqrt{s} - t\right\|_2 + \sqrt{\overline{D}_m \vee \Delta_m}\right\}\right]^q, \end{aligned}$$

*with a constant $C(q)$ depending on $q$ only.*

When $\mu = \mu_s \in \mathcal{Q}_\lambda$, (1.14) becomes

$$(1.15) \quad \mathbb{E}_s\left[H^q(s,\hat{s})\right] \leq C(q)\left[1+\Sigma\right] \inf_{m \in \mathcal{M}} \left\{ \inf_{t \in \overline{S}_m} \left\|\sqrt{s} - t\right\|_2 + \sqrt{\overline{D}_m \vee \Delta_m} \right\}^q.$$

Typical examples for $\mathcal{X}$ and $\lambda$ are $[0,1]^k$ with the Lebesgue measure or $\{1;\ldots;n\}$ with the counting measure. In this last case, the $n$ random variables $\Lambda_{\boldsymbol{X}}(\{i\}) = N_i$ are independent Poisson variables with respective parameters $s_i = s(i)$ and observing $\boldsymbol{X}$ is equivalent to observing a set of $n$ independent Poisson variables with varying parameters, a framework which is usually studied under the name of Poisson regression.

### 1.5. Model selection for Poisson processes, a brief review

Although there have been numerous papers devoted to estimation of the mean measure of a Poisson process, only a few, recently, considered the problem of model selection, the key reference being Reynaud-Bouret [30] with extensions to more general processes in Reynaud-Bouret [31]. A major difference with our approach is her use of the $\mathbb{L}_2(\lambda)$-loss, instead of the Hellinger type loss that we introduce here. It first requires that the unknown mean measure $\mu$ be dominated by $\lambda$ with intensity $s$ and that $s \in \mathbb{L}_2(\lambda)$. Moreover, as we shall show in Section 2.3 the use of the $\mathbb{L}_2$-loss typically requires that $s \in \mathbb{L}_\infty(\lambda)$. This results in rather complicated assumptions but the advantage of this approach is that it is based on penalized projection estimators which can be computed practically while the construction of our estimators is too computationally intensive to be implemented on a computer, as we shall explain below. The same conclusions essentially apply to all other papers dealing with the subject. The approach of Grégoire and Nembé [21], which extends previous results of Barron and Cover [8] about density estimation to that of intensities, has some similarities with ours. The paper by Kolaczyk and Nowak [25] based on penalized maximum likelihood focuses on Poisson regression. Methods which can also be viewed as cases of model selection are those based on the thresholding of the empirical coefficients with respect to some orthonormal basis. It is known that such a procedure is akin to model selection with models spanned by finite subsets of a basis. They have been considered in Kolaczyk [24], Antoniadis, Besbeas and Sapatinas [1], Antoniadis and Sapatinas [2] and Patil and Wood [28].

### 1.6. An overview of the paper

We already justified the introduction of our Hellinger type loss-functions by the properties of likelihood ratio tests and we shall explain, in the next section, why the more popular $\mathbb{L}_2$-risk is not suitable for our purposes, at least if we want to deal with possibly unbounded intensities. To show this, we shall design a general tool for getting lower bounds for intensity estimation, which is merely a version of Assouad's Lemma [3] for Poisson processes. We shall also show that recent results by Rigollet and Tsybakov [32] on aggregation of estimators for density estimation extend straightforwardly to the Poisson case. In Section 3, we briefly recall the general construction of T-estimators introduced in Birgé [9] and apply it to the specific case of Poisson processes. We also provide an illustration based on nonlinear approximating models. Section 4 is devoted to various applications of our method based on families of linear models. This section essentially relies on results from approximation theory about the approximation of different classes of functions (typically smoothness classes) by finite dimensional linear spaces in $\mathbb{L}_2$. We also indicate how to mix different families of models and introduce an asymptotic point of view which allows to consider convergence rates and to make a parallel with density estimation. In Section 5, we deal with aggregation of estimators with some applications to partition selection for histograms. The final Section 6 is devoted to the proof of the most important technical result in this paper, namely the existence and properties of tests between balls of mean measures. This is the key argument which is required to apply the construction of T-estimators to the problem of estimating the mean measure of a Poisson process. It also has other applications, in particular to the study of Bayesian procedures as done, for instance, in Ghosal, Ghosh and van der Vaart [20] and subsequent work of van der Vaart and coauthors.

## 2. Estimation with $\mathbb{L}_2$-loss

### 2.1. From density to intensity estimation

A classical approach to density estimation is based on $\mathbb{L}_2$-loss. We assume that the observations $X_1, \ldots, X_n$ have a density $s_1$ with respect to some dominating measure $\lambda$ and that $s_1$ belongs to the Hilbert space $\mathbb{L}_2(\lambda)$ with scalar product $\langle \cdot, \cdot \rangle$ and norm $\|\cdot\|_2$. Given an estimator $\hat{s}(X_1, \ldots, X_n)$ we define its risk by $\mathbb{E}[\|\hat{s} - s_1\|_2^2]$. In this theory, a central role is played by projection estimators as defined by Cencov [14]. Model selection based on projection estimators has been considered by Birgé and Massart [11]. A more modern treatment can be found in Massart [27]. Thresholding estimators based on wavelet expansions as described in Cohen, DeVore, Kerkyacharian and Picard [15] (see also the many further references therein) can also be viewed as special cases of those. Recently Rigollet and Tsybakov [32] introduced an aggregation method based on projection estimators. Projection estimators have the advantage of simplicity and the drawback or requiring somewhat restrictive assumptions on the density $s_1$ to be estimated, not only that it belongs to $\mathbb{L}_2$ but most of the time to $\mathbb{L}_\infty$. As shown in Birgé [10], Section 5.4.1, the fact that $s_1$ belongs to $\mathbb{L}_\infty$ is essentially a necessary condition to have a control on the $\mathbb{L}_2$-risk of estimators of $s_1$.

As indicated in Baraud and Birgé [4] Section 4.2, there is a parallel between the estimation of a density $s_1$ from $n$ i.i.d. observations and the estimation of the intensity $s = ns_1$ from a Poisson process. This suggests to adapt the known results from density estimation to intensity estimation for Poisson processes. We shall briefly explain how it works, when the Poisson process $\boldsymbol{X}$ has an intensity $s \in \mathbb{L}_\infty(\lambda)$ with $\mathbb{L}_\infty$-norm $\|s\|_\infty$.

The starting point is to observe that, given an element $\varphi \in \mathbb{L}_2(\lambda)$, a natural estimator of $\langle \varphi, s \rangle$ is $\overline{\varphi}(\boldsymbol{X}) = \int \varphi \, d\Lambda_{\boldsymbol{X}} = \sum_{i=1}^N \varphi(X_i)$. It follows from (1.1) that

(2.1) $\quad \mathbb{E}_s[\overline{\varphi}(\boldsymbol{X})] = \langle \varphi, s \rangle \quad \text{and} \quad \text{Var}_s(\overline{\varphi}(\boldsymbol{X})) = \int \varphi^2 s \, d\lambda - \langle \varphi, s \rangle^2 \leq \|s\|_\infty \|\varphi\|_2^2.$

Given a $D$-dimensional linear subspace $S'$ of $\mathbb{L}_2(\lambda)$ with an orthonormal basis $\varphi_1, \ldots, \varphi_D$, we can estimate $s$ by the projection estimator with respect to $S'$:

$$\hat{s}(\boldsymbol{X}) = \sum_{j=1}^D \left[ \sum_{i=1}^N \varphi_j(X_i) \right] \varphi_j.$$

It follows from (2.1) that its risk is bounded by

$$\mathbb{E}_s\left[\|\hat{s}(\boldsymbol{X}) - s\|_2^2\right] \leq \inf_{t \in S'} \|t - s\|_2^2 + \|s\|_\infty D. \tag{2.2}$$

Note that $\hat{s}(\boldsymbol{X})$ is not necessarily an intensity since it may take negative values. This can be fixed: replacing $\hat{s}(\boldsymbol{X})$ by its positive part can only reduce the risk since $s$ is nonnegative.

## 2.2. Aggregation of preliminary estimators

The purpose of this section is to extend some recent results for aggregation of density estimators due to Rigollet and Tsybakov [32] to intensity estimation. The basic tool for aggregation in the context of Poisson processes is the procedure of "thinning" which is the equivalent of sample splitting for i.i.d. observations, see for instance Reiss [29], page 68. Assume that we have at our disposal a Poisson process with mean measure $\mu$: $\Lambda_{\boldsymbol{X}} = \sum_{i=1}^{N} \delta_{X_i}$ and an independent sequence $(Y_i)_{i\geq 1}$ of i.i.d. Bernoulli variables with parameter $p \in (0,1)$. Then the two random measures $\Lambda_{\boldsymbol{X_1}} = \sum_{i=1}^{N} Y_i \delta_{X_i}$ and $\Lambda_{\boldsymbol{X_2}} = \sum_{i=1}^{N} (1-Y_i)\delta_{X_i}$ are two independent Poisson processes with respective mean measures $p\mu$ and $(1-p)\mu$.

Now assume that $\boldsymbol{X}$ is a Poisson process with intensity $s$ with respect to $\lambda$, that $\boldsymbol{X_1}$ and $\boldsymbol{X_2}$ have been derived from $\boldsymbol{X}$ by thinning and that we have at our disposal a finite family $\{\hat{s}_m(\boldsymbol{X_1}), m \in \mathcal{M}\}$ of estimators of $ps$ based on the first process and belonging to $\mathbb{L}_2(\lambda)$. They may be projection estimators or others. These estimators span a $D$-dimensional linear subspace of $\mathbb{L}_2(\lambda)$ with an orthonormal basis $\varphi_1, \ldots, \varphi_D$, $D \leq |\mathcal{M}|$. Working conditionally with respect to $\boldsymbol{X_1}$, we use $\boldsymbol{X_2}$ to build a projection estimator $\tilde{s}(\boldsymbol{X_2})$ of $(1-p)s$ belonging to the linear span of the estimators $\hat{s}_m(\boldsymbol{X_1})$. This is exactly the method used by Rigollet and Tsybakov [32] for density estimation and the proof of their Theorem 2.1 extends straightforwardly to Poisson processes to give

**Theorem 2.** *The aggregated estimator $\tilde{s}$ based on the processes $\boldsymbol{X_1}$ and $\boldsymbol{X_2}$ by thinning of $\boldsymbol{X}$ satisfies*

$$\mathbb{E}_s\left[\|\tilde{s}(\boldsymbol{X}) - (1-p)s\|_2^2\right] \leq \mathbb{E}_s\left[\inf_{\theta \in \mathbb{R}^{\mathcal{M}}} \left\|ps - \sum_{m \in \mathcal{M}} \theta_m \hat{s}_m(\boldsymbol{X_1})\right\|_2^2\right] + (1-p)\|s\|_\infty |\mathcal{M}|. \tag{2.3}$$

Setting $\hat{s}(\boldsymbol{X}) = \tilde{s}(\boldsymbol{X})/(1-p)$ leads to

$$\mathbb{E}_s\left[\|\hat{s}(\boldsymbol{X}) - s\|_2^2\right] \leq \frac{1}{(1-p)^2} \inf_{m \in \mathcal{M}} \mathbb{E}_s\left[\|ps - \hat{s}_m(\boldsymbol{X_1})\|_2^2\right] + \frac{\|s\|_\infty |\mathcal{M}|}{1-p}.$$

If we start with a finite family $\{\overline{S}_m, m \in \mathcal{M}\}$ of finite-dimensional linear subspaces of $\mathbb{L}_2(\lambda)$ with respective dimensions $D_m$, we may choose for $\hat{s}_m(\boldsymbol{X_1})$ the projection estimator based on $\overline{S}_m$ with risk bounded by (2.2)

$$\mathbb{E}_s\left[\|\hat{s}_m(\boldsymbol{X_1}) - ps\|_2^2\right] \leq \inf_{t \in \overline{S}_m} \|t - ps\|_2^2 + p\|s\|_\infty D_m = p^2 \inf_{t \in \overline{S}_m} \|t - s\|_2^2 + p\|s\|_\infty D_m.$$

Choosing $p = 1/2$, we conclude that

$$\mathbb{E}_s\left[\|\hat{s}(\boldsymbol{X}) - s\|_2^2\right] \leq \inf_{m \in \mathcal{M}} \left\{\inf_{t \in \overline{S}_m} \|t - s\|_2^2 + 2\|s\|_\infty D_m\right\} + 2\|s\|_\infty |\mathcal{M}|.$$

## *2.3. Lower bounds for intensity estimation*

It is rather inconvenient to get risk bounds involving the unknown and possibly very large $\mathbb{L}_\infty$-norm of $s$ and this problem becomes even more serious if $s$ does not belong to $\mathbb{L}_\infty(\lambda)$. It is, unfortunately, impossible to avoid this problem when dealing with the $\mathbb{L}_2$-loss. To show this, let us start with a version of Assouad's Lemma [3] for Poisson processes.

**Lemma 1.** *Let $\mathcal{S}_D = \{s_\delta, \delta \in \mathcal{D}\} \subset \mathbb{L}_1^+(\lambda)$ be a family of intensities indexed by $\mathcal{D} = \{0;1\}^D$ and $\Delta$ be the Hamming distance on $\mathcal{D}$ given by $\Delta(\delta, \delta') = \sum_{j=1}^D |\delta_j - \delta_j'|$. Let $\mathcal{C}$ be the subset of $\mathcal{D} \times \mathcal{D}$ defined by*

$$\mathcal{C} = \{(\delta, \delta') \,|\, \exists k, 1 \leq k \leq D \quad \text{with } \delta_k = 0, \delta_k' = 1 \quad \text{and} \quad \delta_j = \delta_j' \text{ for } j \neq k\}.$$

*Then for any estimator $\hat{\delta}(\boldsymbol{X})$ with values in $\mathcal{D}$,*

$$(2.4) \qquad \sup_{\delta \in \mathcal{D}} \mathbb{E}_{s_\delta}\left[\Delta\left(\hat{\delta}(\boldsymbol{X}), \delta\right)\right] \geq \frac{D}{4}\left(\frac{1}{|\mathcal{C}|} \sum_{(\delta, \delta') \in \mathcal{C}} \exp\left[-2H^2(s_\delta, s_{\delta'})\right]\right).$$

*If, moreover, $\mathcal{S}_D \subset \mathbb{L} \subset \mathbb{L}_1^+(\lambda)$ and $\mathbb{L}$ is endowed with a metric $d$ satisfying $d^2(s_\delta, s_{\delta'}) \geq \theta \Delta(\delta, \delta')$ for all $\delta, \delta' \in \mathcal{D}$ and some $\theta > 0$, then for any estimator $\hat{s}(\boldsymbol{X})$ with values in $\mathbb{L}$,*

$$(2.5) \qquad \sup_{s \in \mathcal{S}_D} \mathbb{E}_s\left[d^2\left(\hat{s}(\boldsymbol{X}), s\right)\right] \geq \frac{D\theta}{16}\left(\frac{1}{|\mathcal{C}|} \sum_{(\delta, \delta') \in \mathcal{C}} \exp\left[-2H^2(s_\delta, s_{\delta'})\right]\right).$$

*Proof.* To get (2.4) it suffices to find a lower bound for

$$R_B = 2^{-D} \sum_{\delta \in \mathcal{D}} \mathbb{E}_{s_\delta}\left[\Delta\left(\hat{\delta}, \delta\right)\right] = 2^{-D} \sum_{\delta \in \mathcal{D}} \int \sum_{k=1}^D \left|\hat{\delta}_k - \delta_k\right| dQ_{s_\delta},$$

since the left-hand side of (2.4) is at least as large as the average risk $R_B$. It follows from the proof of Lemma 2 in Birgé [10] with $n = 1$ that

$$R_B \geq 2^{-D} \sum_{(\delta, \delta') \in \mathcal{C}} \left[1 - \sqrt{1 - \rho^2(Q_{s_\delta}, Q_{s_{\delta'}})}\right] \geq 2^{-D-1} \sum_{(\delta, \delta') \in \mathcal{C}} \rho^2(Q_{s_\delta}, Q_{s_{\delta'}}).$$

Then (2.4) follows from (1.6) since $|\mathcal{C}| = D2^{D-1}$. Let now $\hat{s}(\boldsymbol{X})$ be an estimator with values in $\mathbb{L}$ and set $\hat{\delta}(\boldsymbol{X}) \in \mathcal{D}$ to satisfy $d(\hat{s}, s_{\hat{\delta}}) = \inf_{\delta \in \mathcal{D}} d(\hat{s}, s_\delta)$ so that, whatever $\delta \in \mathcal{D}$, $d(s_{\hat{\delta}}, s_\delta) \leq 2d(\hat{s}, s_\delta)$. It then follows from our assumptions that

$$\sup_{\delta \in \mathcal{D}} \mathbb{E}_{s_\delta}\left[d^2\left(\hat{s}, s_\delta\right)\right] \geq \frac{1}{4} \sup_{\delta \in \mathcal{D}} \mathbb{E}_{s_\delta}\left[d^2\left(s_{\hat{\delta}}, s_\delta\right)\right] \geq \frac{\theta}{4} \sup_{\delta \in \mathcal{D}} \mathbb{E}_{s_\delta}\left[\Delta\left(\hat{\delta}(\boldsymbol{X}), \delta\right)\right]$$

and (2.5) follows from (2.4). □

The simplest application of this lemma corresponds to the case $D = 1$ which, in its simplest form, dates back to Le Cam [26]. We consider only two intensities $s_0$ and $s_1$ so that $\theta = d^2(s_0, s_1)$ and (2.5) gives, whatever the estimator $\hat{s}(\boldsymbol{X})$,

$$(2.6) \qquad \max_{i=0,1} \mathbb{E}_{s_i}\left[d^2\left(\hat{s}(\boldsymbol{X}), s_i\right)\right] \geq \frac{d^2(s_0, s_1)}{16} \exp\left[-2H^2(s_0, s_1)\right].$$

Another typical application of the previous lemma to intensities on $[0, 1]$ uses the following construction of a suitable set $\mathcal{S}_D$.

**Lemma 2.** *Let $D$ be a positive integer and $g$ be a function on $\mathbb{R}$ with support on $[0, D^{-1})$ satisfying*

$$0 \leq g(x) \leq 1 \quad \text{for all } x \quad \text{and} \quad \int_0^{D^{-1}} g^2(x)\, dx = a > 0.$$

*Set, for $1 \leq j \leq D$ and $0 \leq x \leq 1$, $g_j(x) = g(x - D^{-1}(j-1))$ and, for $\delta \in \mathcal{D}$, $s_\delta(x) = a^{-1}[1 + \sum_{j=1}^D (\delta_j - 1/2)g_j(x)]$. Then $\|s_\delta - s_{\delta'}\|_2^2 = a^{-1}\Delta(\delta, \delta')$ and $H^2(s_\delta, s_{\delta'}) \geq \Delta(\delta, \delta')/8$ for all $\delta, \delta' \in \mathcal{D}$. Moreover,*

$$(2.7) \quad |\mathcal{C}|^{-1} \sum_{(\delta,\delta') \in \mathcal{C}} \exp\left[-2H^2(s_\delta, s_{\delta'})\right] \geq \exp[-2/7].$$

*Proof.* The first equality is clear. Let us then observe that our assumptions on $g$ imply that $1 - g^2(x)/7 \leq \sqrt{1 - g^2(x)/4} \leq 1 - g^2(x)/8$, hence, since the functions $g_j$ have disjoint supports and are translates of $g$,

$$H^2(s_\delta, s_{\delta'}) = (2a)^{-1} \sum_{j=1}^D |\delta_j - \delta'_j| \int_0^{D^{-1}} \left[\sqrt{1 + g(x)/2} - \sqrt{1 - g(x)/2}\right]^2 dx$$

$$= a^{-1} \sum_{j=1}^D |\delta_j - \delta'_j| \int_0^{D^{-1}} \left[1 - \sqrt{1 - g^2(x)/4}\right] dx = c\Delta(\delta, \delta'),$$

with $1/8 \leq c \leq 1/7$. The conclusions follow. □

**Corollary 1.** *For each positive integer $D$ and $L \geq 3D/2$, one can find a finite set $\mathcal{S}_D$ of intensities with the following properties:*
  *(i) it is a subset of some $D$-dimensional affine subspace of $\mathbb{L}_2([0,1], dx)$;*
  *(ii) $\sup_{s \in \mathcal{S}_D} \|s\|_\infty \leq L$;*
  *(iii) for any estimator $\hat{s}(\boldsymbol{X})$ with values in $\mathbb{L}_2([0,1], dx)$ based on a Poisson process $\boldsymbol{X}$ with intensity $s$,*

$$(2.8) \quad \sup_{s \in \mathcal{S}_D} \mathbb{E}_s\left[\|\hat{s} - s\|_2^2\right] \geq (DL/24)\exp[-2/7].$$

*Proof.* Let us set $\theta = 2L/3 \geq D$ and apply the construction of Lemma 2 with $g(x) = \sqrt{D/\theta}\, \mathbb{1}_{[0,1/D)}$, hence $a = \theta^{-1}$. This results in the set $\mathcal{S}_D$ with $\|s_\delta\|_\infty \leq \theta\left[1 + (1/2)\sqrt{D/\theta}\right] \leq 3\theta/2 = L$ for all $\delta \in \mathcal{D}$ as required. Moreover $\|s_\delta - s_{\delta'}\|_2^2 = \theta\Delta(\delta, \delta')$. Then we use Lemma 1 with $d$ being the distance corresponding to the norm in $\mathbb{L}_2([0,1], dx)$ and (2.5) together with (2.7) result in (2.8). □

This result implies that, if we want to use the squared $\mathbb{L}_2$-norm as a loss function, whatever the choice of our estimator there is no hope to find risk bounds that are independent of the $\mathbb{L}_\infty$-norm of the underlying intensity, even if this intensity belongs to a finite-dimensional affine space. This provides an additional motivation for the introduction of loss functions based on the distance $H$.

## 3. T-estimators for Poisson processes

### 3.1. Some notations

Throughout this paper, we observe a Poisson process $\boldsymbol{X}$ on $\mathcal{X}$ with unknown mean measure $\mu$ belonging to the metric space $(\mathcal{Q}_+(\mathcal{X}), H)$ and have at hand some reference measure $\lambda$ on $\mathcal{X}$ so that $\mu = \mu_s + \mu^\perp$ with $\mu_s \in \mathcal{Q}_\lambda$, $s \in \mathbb{L}_1^+(\lambda)$ and $\mu^\perp$ orthogonal to $\lambda$. We denote by $\|\cdot\|_i$ the norm in $\mathbb{L}_i(\lambda)$ for $1 \leq i \leq \infty$ and by $d_2$ the distance corresponding to the norm $\|\cdot\|_2$. We always denote by $s$ the intensity of the part of $\mu$ which is dominated by $\lambda$ and set $s_1 = s/\mu_s(\mathcal{X})$. We also systematically identify $\mathcal{Q}_\lambda$ with $\mathbb{L}_1^+(\lambda)$ via the mapping $t \mapsto \mu_t$, writing $t$ as a shorthand for $\mu_t \in \mathcal{Q}_\lambda$. We write $H(s, S')$ for $\inf_{t \in S'} H(s, t)$, $a \vee b$ and $a \wedge b$ for the maximum and the minimum respectively of $a$ and $b$, $|A|$ for the cardinality of a finite set $A$ and $\mathbb{N}^\star = \mathbb{N} \setminus \{0\}$ for the set of positive integers. In the sequel $C$ (or $C', C_1, \ldots$) denote constants that may vary from line to line, the form $C(a, b)$ meaning that $C$ is not a universal constant but depends on some parameters $a$ and $b$.

### 3.2. Definition and properties of T-estimators

In order to explain our method of estimation and model selection, we need to recall some general results from Birgé [9] about T-estimators that we shall specialize to the specific framework of this paper. Let $(M, d)$ be some metric space and $\mathcal{B}(t, r)$ denote the open ball of center $t$ and radius $r$ in $M$.

**Definition 1.** *A subset $S'$ of the metric space $(M, d)$ is called* a D-model *with parameters $\eta, D$ and $B'$ ($\eta, B', D > 0$) if*

$$(3.1) \quad |S' \cap \mathcal{B}(t, x\eta)| \leq B' \exp\left[Dx^2\right] \quad \text{for all } x \geq 2 \text{ and } t \in M.$$

Note that this implies that $S'$ is at most countable.

To estimate the unknown mean measure $\mu$ of the Poisson process $\boldsymbol{X}$, we introduce a finite or countable family $\{S_m, m \in \mathcal{M}\}$ of D-models in $(\mathcal{Q}_\lambda, H)$ with respective parameters $\eta_m, D_m$ and $B'$ and assume that

$$(3.2) \quad \text{for all } m \in \mathcal{M}, \quad D_m \geq 1/2 \quad \text{and} \quad \eta_m^2 \geq (84 D_m)/5,$$

and

$$(3.3) \quad \sum_{m \in \mathcal{M}} \exp\left[-\eta_m^2/84\right] = \Sigma < +\infty.$$

Then we set $S = \bigcup_{m \in \mathcal{M}} S_m$ and, for each $t \in S$,

$$(3.4) \quad \eta(t) = \inf\{\eta_m \,|\, m \in \mathcal{M} \text{ and } S_m \ni t\}.$$

**Remark.** Note that if we choose for $\{S_m, m \in \mathcal{M}\}$ a family of D-models in $(\mathcal{Q}_+(\mathcal{X}), H)$, $S$ is countable and therefore dominated by some measure $\lambda$ that we can always take as our reference measure. This gives an a posteriori justification for the choice of a family of models $S_m \subset \mathcal{Q}_\lambda$.

Given two distinct points $t, u \in \mathcal{Q}_\lambda$ we define a test function $\psi(\boldsymbol{X})$ between $t$ and $u$ as a measurable function from $\mathcal{X}$ to $\{t, u\}$, $\psi(\boldsymbol{X}) = t$ meaning deciding $t$ and $\psi(\boldsymbol{X}) = u$ meaning deciding $u$. In order to define a T-estimator, we need a family of test functions $\psi_{t,u}(\boldsymbol{X})$ between distinct points $t, u \in S$ with some special properties. The following proposition, to be proved in Section 6 warrants their existence.

**Proposition 1.** *Given two distinct points $t, u \in S$ there exists a test $\psi_{t,u}$ between $t$ and $u$ which satisfies*

$$\sup_{\{\mu \in \mathcal{Q}_+(\mathcal{X}) \,|\, H(\mu, \mu_t) \leq H(t,u)/4\}} \mathbb{P}_\mu[\psi_{t,u}(\boldsymbol{X}) = u]$$
$$\leq \exp\left[-\left(H^2(t,u) - \eta^2(t) + \eta^2(u)\right)/4\right],$$

$$\sup_{\{\mu \in \mathcal{Q}_+(\mathcal{X}) \,|\, H(\mu, \mu_u) \leq H(t,u)/4\}} \mathbb{P}_\mu[\psi_{t,u}(\boldsymbol{X}) = t]$$
$$\leq \exp\left[-\left(H^2(t,u) - \eta^2(u) + \eta^2(t)\right)/4\right],$$

and for all $\mu \in \mathcal{Q}_+(\mathcal{X})$,

(3.5) $\quad \mathbb{P}_\mu[\psi_{t,u}(\boldsymbol{X}) = u] \leq \exp\left[\left(16H^2(\mu, \mu_t) + \eta^2(t) - \eta^2(u)\right)/4\right].$

To build a T-estimator, we proceed as follows. We consider a family of tests $\psi_{t,u}$ indexed by the two-points subsets $\{t, u\}$ of $S$ with $t \neq u$ that satisfy the conclusions of Proposition 1 and we set $\mathcal{R}_t = \{u \in S, u \neq t \,|\, \psi_{t,u}(\boldsymbol{X}) = u\}$ for each $t \in S$. Then we define the random function $\mathcal{D}_{\boldsymbol{X}}$ on $S$ by

$$\mathcal{D}_{\boldsymbol{X}}(t) = \begin{cases} \sup_{u \in \mathcal{R}_t} \{H(t,u)\} & \text{if } \mathcal{R}_t \neq \emptyset; \\ 0 & \text{if } \mathcal{R}_t = \emptyset. \end{cases}$$

We call *T-estimator* derived from $S$ and the family of tests $\psi_{t,u}(\boldsymbol{X})$ any measurable minimizer of the function $t \mapsto \mathcal{D}_{\boldsymbol{X}}(t)$ from $S$ to $[0, +\infty]$ so that $\mathcal{D}_{\boldsymbol{X}}(\hat{s}(\boldsymbol{X})) = \inf_{t \in S} \mathcal{D}_{\boldsymbol{X}}(t)$. Such a minimizer need not exist in general but it actually exists under our assumptions.

**Theorem 3.** *Let $S = \bigcup_{m \in \mathcal{M}} S_m \subset \mathcal{Q}_\lambda$ be a finite or countable family of D-models in $(\mathcal{Q}_\lambda, H)$ with respective parameters $\eta_m, D_m$ and $B'$ satisfying (3.2) and (3.3). Let $\{\psi_{t,u}\}$ be a family of tests indexed by the two-points subsets $\{t, u\}$ of $S$ with $t \neq u$ and satisfying the conclusions of Proposition 1. Whatever $\mu \in \mathcal{Q}_+(\mathcal{X})$, $\mathbb{P}_\mu$-a.s. there exists at least one T-estimator $\hat{s} = \hat{s}(\boldsymbol{X}) \in S$ derived fom this family of tests and any of them satisfies, for all $s' \in S$,*

(3.6) $\quad \mathbb{P}_\mu\left[H(s', \hat{s}) > y\right] < (B'\Sigma/7) \exp\left[-y^2/6\right] \quad \text{for } y \geq 4[H(\mu, \mu_{s'}) \vee \eta(s')].$

*Setting $\hat{\mu}(\boldsymbol{X}) = \hat{s}(\boldsymbol{X}) \cdot \lambda$ and $\mu = \mu_s + \mu^\perp$ with $\mu_s \in \mathcal{Q}_\lambda$ and $\mu^\perp$ orthogonal to $\lambda$, we also get*

(3.7) $\quad \mathbb{E}_\mu\left[H^q\left(\mu, \hat{\mu}(\boldsymbol{X})\right)\right] \leq C(q)[1 + B'\Sigma] \inf_{m \in \mathcal{M}} \left\{H(s, S_m) + \eta_m + \sqrt{\mu^\perp(\mathcal{X})}\right\}^q$

*and, for intensity estimation when $\mu = \mu_s$,*

(3.8) $\quad \mathbb{E}_s\left[H^q\left(s, \hat{s}(\boldsymbol{X})\right)\right] \leq C(q)[1 + B'\Sigma] \inf_{m \in \mathcal{M}} \{H(s, S_m) + \eta_m\}^q.$

*Proof.* It follows from Theorem 5 in Birgé [9] with $a = 1/4$, $B = 1$, $\kappa = 4$ and $\kappa' = 16$ that T-estimators do exist, satisfy (3.6) and have a risk which is bounded, for $q \geq 1$, by

(3.9) $\quad \mathbb{E}_\mu\left[H^q\left(\mu, \hat{\mu}(\boldsymbol{X})\right)\right] \leq C(q)[1 + B'\Sigma] \inf_{m \in \mathcal{M}} \left\{\left(\inf_{t \in S_m} H(\mu, \mu_t)\right) \vee \eta_m\right\}^q.$

In Birgé [9], the proof of the existence of T-estimators when $\mathcal{M}$ is infinite was given only for the case that the tests $\psi_{t,u}(\boldsymbol{X})$ have a special form, namely $\psi_{t,u}(\boldsymbol{X}) = u$ when $\gamma(u, \boldsymbol{X}) < \gamma(t, \boldsymbol{X})$ and $\psi_{t,u}(\boldsymbol{X}) = t$ when $\gamma(u, \boldsymbol{X}) > \gamma(t, \boldsymbol{X})$ for some suitable function $\gamma$. A minor modification of the proof extends the result to the general situation based on the assumption that (3.5) holds. It is indeed enough to use (3.5) to modify the proof of (7.18) of Birgé [9] in order to get instead

$$\mathbb{P}_\mu \left[ \exists t \in S \text{ with } \psi_{s',t}(\boldsymbol{X}) = 1 \text{ and } \eta(t) \geq y \right] \xrightarrow[y \to +\infty]{} 0.$$

The existence of $\hat{s}(\boldsymbol{X})$ then follows straightforwardly. Since $H^2(\mu, \mu_t) = H^2(s,t) + \mu^\perp(\mathcal{X})/2$, (3.7) follows from (3.9). □

It follows from (3.7) that the problem of estimating $\mu$ with T-estimators always reduces to intensity estimation once a reference measure $\lambda$ has been chosen. A comparison of the risk bounds (3.7) and (3.8) shows that the performance of the estimator $\hat{s}(\boldsymbol{X})$ is connected to the choice of the models in $\mathbb{L}_1^+(\lambda)$, the component $\mu^\perp(\mathcal{X})$ of the risk depending only on $\lambda$. We might as well assume that $\mu^\perp(\mathcal{X})$ is known since this would not change anything concerning the performance of the T-estimators for a given $\lambda$. This is why we shall essentially focus, in the sequel, on intensity estimation.

### 3.3. An application to multivariate intensities

Let us first illustrate Theorem 3 by an application to the estimation of the unknown intensity $s$ (with respect to the Lebesgue measure $\lambda$) of a Poisson process on $\mathcal{X} = [-1, 1]^k$. For this, we introduce a family of non-linear models related to neural nets which were popularized in the 90's by Barron [5, 6] and other authors in view of their nice approximation properties with respect to functions of several variables. These models have already been studied in detail in Sections 3.2.2 and 4.2.2 of Barron, Birgé and Massart [7] and we shall therefore refer to this paper for their properties. We start with a family of functions $\phi_w(x) \in \mathbb{L}_\infty([-1,1]^k)$ indexed by a parameter $w$ belonging to $\mathbb{R}^{k'}$ and satisfying

(3.10) $$|\phi_w(x) - \phi_{w'}(x)| \leq |w - w'|_1 \quad \text{for all } x \in [-1,1]^k,$$

where $|\cdot|_1$ denotes the $l^1$-norm on $\mathbb{R}^{k'}$. Various examples of such families are given in Barron, Birgé and Massart [7] and one can, for instance, set $\phi_w(x) = \psi(a'x - b)$ with $\psi$ a univariate Lipschitz function, $a \in \mathbb{R}^k, b \in \mathbb{R}$ and $w = (a,b) \in \mathbb{R}^{k+1}$.

We set $\mathcal{M} = (\mathbb{N} \setminus \{0,1\})^3$ and for $m = (J, R, B) \in \mathcal{M}$ we consider the subset of $\mathbb{L}_\infty([-1,1]^k)$ defined by

$$S'_m = \left\{ \sum_{j=1}^J \beta_j \phi_{w_j}(x) \ \bigg| \ \sum_{j=1}^J |\beta_j| \leq R \text{ and } |w_j|_1 \leq B \text{ for } 1 \leq j \leq J \right\}.$$

As shown in Lemma 5 of Barron, Birgé and Massart [7], such a model can be approximated by a finite subset $T_m$. More precisely, one can find a subset $T_m$ of $S'_m$ with cardinality bounded by $[2e(2RB+1)]^{J(k'+1)}$ and such that if $u \in S'_m$, there exists some $t \in T_m$ such that $\|t - u\|_\infty \leq 1$. Defining $S_m$ as $\{t^2, t \in T_m\}$, we get the following property:

**Lemma 3.** *For $m = (J, R, B) \in (\mathbb{N} \setminus \{0,1\})^3$, we set $\eta_m^2 = 42J(k'+1)\log(RB)$. Then $S_m$ is a D-model with parameters $\eta_m, D_m = [J(k'+1)/4]\log[2e(2RB+1)]$ and 1 in the metric space $(\mathbb{L}_1^+(\lambda), H)$ and (3.2) and (3.3) are satisfied. Moreover, for any $s \in \mathbb{L}_1^+(\lambda)$,*

(3.11) $$\sqrt{2}H(s, S_m) \leq \inf_{t \in S'_m} \left\| \sqrt{s} - t \right\|_2 + 2^{k/2}.$$

*Proof.* Since $|S_m| \leq |T_m|$, to show that $S_m$ is a D-model with the given parameters it is enough to prove, in view of (3.1), that $|T_m| \leq \exp[4D_m]$, which is clear. That $\eta_m^2/84 \geq D_m/5$ follows from $\log[2e(2RB+1)] \leq 4\log(RB)$ since $RB \geq 4$. Moreover, since $k'+1 \geq 2$, $\eta_m^2 \geq 84J\log(RB)$, hence

$$\sum_{m \in \mathcal{M}} \exp\left[-\frac{\eta_m^2}{84}\right] \leq \sum_{J \geq 2}\left(\sum_{n \geq 2} n^{-J}\right)^2 \leq \sum_{J \geq 2}\left(\int_{3/2}^{+\infty} x^{-J}\,dx\right)^2,$$

so that (3.3) holds. Let now $u \in S'_m$. There exists $t \in T_m$ such that $\|t-u\|_\infty \leq 1$, hence $\|\sqrt{s}-t\|_2 \leq \|\sqrt{s}-u\|_2 + 2^{k/2}$. Then $t^2 \in S_m$ and since $\|\sqrt{s}-\sqrt{t^2}\|_2 \leq \|\sqrt{s}-t\|_2$, (3.11) follows. □

Let now $\hat{s}(\boldsymbol{X})$ be a T-estimator derived from the family of D-models $\{S_m, m \in \mathcal{M}\}$. By Theorem 3 and Lemma 3, it satisfies

$$\mathbb{E}_s\left[H^2(s, \hat{s}(\boldsymbol{X}))\right] \leq C \inf_{m \in \mathcal{M}}\left\{\inf_{t \in S'_m}\|\sqrt{s}-t\|_2^2 + 2^k + \eta_m^2\right\}$$

(3.12)
$$\leq C(k,k') \inf_{m \in \mathcal{M}}\left\{\inf_{t \in S'_m}\|\sqrt{s}-t\|_2^2 + J\log(RB)\right\}.$$

The approximation properties of the models $S'_m$ with respect to different classes of functions have been described in Barron, Birgé and Massart [7]. They allow to bound $\inf_{t \in S'_m}\|\sqrt{s}-t\|_2$ when $\sqrt{s}$ belongs to such classes so that corresponding risk bounds can be derived from (3.12).

### 3.4. Model selection based on linear models

#### 3.4.1. Deriving D-models from linear spaces

In order to apply Theorem 3 we need to introduce suitable families of D-models $S_m$ in $(\mathcal{Q}_\lambda, H)$ with good approximation properties with respect to the unknown $s$. More precisely, it follows from (3.7) and (1.10) that they should provide approximations of $\sqrt{s}$ in $\mathbb{L}_2^+(\lambda)$. Good approximating sets for elements of $\mathbb{L}_2^+(\lambda)$ are provided by approximation theory and some recipes to derive D-models from such sets have been given in Section 6 of Birgé [9]. Most results about approximation of functions in $\mathbb{L}_2(\lambda)$ deal with finite dimensional linear spaces or unions of such spaces and their approximation properties with respect to different classes (typically smoothness classes) of functions. We therefore focus here on such linear subspaces of $\mathbb{L}_2(\lambda)$. To translate their properties in terms of D-models, we shall invoke the following proposition.

**Proposition 2.** *Let $\overline{S}$ be a $k$-dimensional linear subspace of $\mathbb{L}_2(\lambda)$ and $\delta > 0$. One can find a subset $S'$ of $\mathcal{Q}_\lambda$ which is a D-model in the metric space $(\mathcal{Q}_\lambda, H)$ with parameters $\delta, 9k$ and $1$ and such that, for any intensity $s \in \mathbb{L}_1^+(\lambda)$,*

$$H(s, S') \leq 2.2\left[\inf_{t \in \overline{S}}\|\sqrt{s}-t\|_2 + \delta\right].$$

*Proof.* Let us denote by $\mathcal{B}_H$ and $\mathcal{B}_2$ the open balls in the metric spaces $(\mathbb{L}_1^+(\lambda), H)$ and $(\mathbb{L}_2(\lambda), d_2)$ respectively. It follows from Proposition 8 of Birgé [9] that one can find a subset $T$ of $\overline{S}$ which is a D-model of $(\mathbb{L}_2(\lambda), d_2)$ with parameters $\delta, k/2$ and 1 and such that, whatever $u \in \mathbb{L}_2(\lambda)$, $d_2(u,T) \leq d_2(u,\overline{S}) + \delta$. It follows that

$$(3.13) \quad \left|T \cap \mathcal{B}_2\left(t, 3r'\sqrt{2}\right)\right| \leq \exp\left[9k(r'/\delta)^2\right] \quad \text{for } r' \geq 2\delta \quad \text{and } t \in \mathbb{L}_2(\lambda).$$

Moreover, if $t \in T$, $\overline{\pi}(t) = \max\{t,0\}$ belongs to $\mathbb{L}_2^+(\lambda)$ and satisfies $d_2(u, \overline{\pi}(t)) \leq d_2(u,t)$ for any $u \in \mathbb{L}_2^+(\lambda)$. We may therefore apply Proposition 12 of Birgé [9] with $(M',d) = (\mathbb{L}_2(\lambda), d_2)$, $M_0 = \mathbb{L}_2^+(\lambda)$, $\lambda = 1$, $\varepsilon = 1/10$, $\eta = 4\sqrt{2}\delta$ and $r = r'\sqrt{2}$ to get a subset $S$ of $\overline{\pi}(T) \subset \mathbb{L}_2^+(\lambda)$ such that

$$(3.14) \quad |S \cap \mathcal{B}_2\left(t, r'\sqrt{2}\right)| \leq |T \cap \mathcal{B}_2\left(t, 3r'\sqrt{2}\right)| \vee 1 \quad \text{for all } t \in \mathbb{L}_2(\lambda) \text{ and } r' \geq 2\delta$$

and $d_2(u,S) \leq 3.1 d_2(u,T)$ for all $u \in \mathbb{L}_2^+(\lambda)$. Setting $S' = \{t^2 \cdot \lambda, t \in S\} \subset \mathcal{Q}_\lambda$ and using (1.10), we deduce from (3.13) and (3.14) that

$$|S' \cap \mathcal{B}_H(\mu_t, r')| \leq \exp\left[9k(r'/\delta)^2\right] \quad \text{for } r' \geq 2\delta \quad \text{and } \mu_t \in \mathcal{Q}_\lambda,$$

hence $S'$ is a D-model in $(\mathcal{Q}_\lambda, H)$ with parameters $\delta, 9k$ and 1, and

$$H(s, S') \leq \left(3.1/\sqrt{2}\right) d_2\left(\sqrt{s}, T\right) < 2.2\left[d_2\left(\sqrt{s}, \overline{S}\right) + \delta\right]. \qquad \square$$

We are now in a position to prove Theorem 1. For each $m$, let us fix $\eta_m^2 = 84[\Delta_m \vee (9\overline{D}_m/5)]$ and use Proposition 2 to derive from $\overline{S}_m$ a D-model $S_m$ with parameters $\eta_m, D_m = 9\overline{D}_m$ and 1 which also satisfies

$$H(s, S_m) \leq 2.2 \left[\inf_{t \in \overline{S}_m} \left\|\sqrt{s} - t\right\|_2 + \eta_m\right].$$

It follows from the definition of $\eta_m$ that (3.2) and (3.3) are satisfied so that Theorem 3 applies. The conclusion immediately follows from (3.7).

### 3.4.2. *About the computation of T-estimators*

We already mentioned that the relevance of T-estimators is mainly of a theoretical nature because of the difficulty of their implementation. Let us give here a simple illustrative example based on a single linear approximating space $\overline{S}$ for $\sqrt{s}$, of dimension $k$. To try to get a practical implementation, we shall use a simple discretization strategy. The first step is to replace $\overline{S}$, that we identify to $\mathbb{R}^k$ via the choice of a basis, by $\theta \mathbb{Z}^k$. This provides an $\eta$-net for $\mathbb{R}^k$ with respect to the the Euclidean distance, with $\eta^2 = k(\theta/2)^2$. Let us concentrate here on the case of a large value of $\Gamma^2 = \int s\, d\lambda$ in order to have a large number of observations since $N$ has a Poisson distribution with parameter $\Gamma^2$. In particular, we shall asume that $\Gamma^2$ (which plays the role of the number of observations as we shall see in Section 4.6) is much larger than $k$. It is useless, in such a case, to use the whole of $\theta \mathbb{Z}^k$ to approximate $\sqrt{s}$ since the closest point to $\sqrt{s}$ belongs to $\mathcal{B}(0, \Gamma + \eta)$. Of course, $\Gamma$ is unknown, but when it is large it can be safely estimated by $\sqrt{N}$ in view of the concentration properties of Poisson variables. Let us therefore assume that $N \geq \Gamma^2/2 \geq 2k$. A reasonable approximating set for $\sqrt{s}$ is therefore $T = \mathcal{B}(0, \sqrt{2N} + \eta) \cap \theta \mathbb{Z}^k$ and since our final model $S$ should be a subset of $\mathbb{L}_2^+(\lambda)$, we can take $S = \{t \vee 0, t \in T\}$ so that $d_2(\sqrt{s}, S) \leq d_2(\sqrt{s}, T) \leq d_2(\sqrt{s}, \overline{S} + \eta)$. It follows from Lemma 5 of Birgé [9] that

$$|S| \leq |T| \leq \frac{(\pi e/2)^{k/2}}{\sqrt{\pi k}} \left(\frac{2\sqrt{2N} + 2\eta}{\theta \sqrt{k}} + 1\right)^k < K = \left[c\left(\sqrt{2N}\eta^{-1} + 1\right)\right]^k,$$

with $c = \sqrt{\pi e/2} \sim 2.07$. This implies that $S$ is a D-model with parameters $\eta, (\log K)/4$ and 1. In order that (3.2) be satisfied, we need that $\eta^2 \geq 4.2 \log K$. If we choose $\eta^2 = 4.2k \log(c(\sqrt{N/k} + 1))$, this inequality holds since $\eta \geq 2\sqrt{k}$, hence $K \leq [c(\sqrt{N/k} + 1)]^k$. The number of tests required for building the T-estimator is $|S|(|S| - 1) < K^2$. For $N$ of the order of 100 and $k$ as small as 5, $K^2$ is of the order of $10^{10}$. This toy example illustrates the difficulty of implementing the algorithm. More realistic ones would be much worse.

## 4. Applications with linear models

We now assume that $\mu = \mu_s = s \cdot \lambda$ and focus on the estimation of the intensity $s$ by model selection, starting with linear models in $\mathbb{L}_2(\lambda)$ that possess good approximating properties with respect to $\sqrt{s}$.

### 4.1. Adaptation in Besov spaces

It is now well-known that wavelet bases are very good tools for representing smooth functions in $\mathbb{L}_2([0,1]^l, dx)$. In particular, given a suitable wavelet basis $\{\varphi_{j,\boldsymbol{k}}, j \geq -1, \boldsymbol{k} \in \Lambda(j)\}$ with $|\Lambda(-1)| \leq \Gamma$ and $2^{jl} \leq |\Lambda(j)| \leq \Gamma 2^{jl}$ for all $j \geq 0$ any function $f \in \mathbb{L}_2([0,1]^l, dx)$ can be written as $f = \sum_{j=-1}^{\infty} \sum_{\boldsymbol{k} \in \Lambda(j)} \beta_{j,\boldsymbol{k}} \varphi_{j,\boldsymbol{k}}$. Moreover $f$ belongs to the Besov space $B_{p,\infty}^\alpha([0,1]^l)$ if and only if

$$(4.1) \quad \sup_{j \geq 0} 2^{j(\alpha + \frac{l}{2} - \frac{l}{p})} \left( \sum_{\boldsymbol{k} \in \Lambda(j)} |\beta_{j,\boldsymbol{k}}|^p \right)^{\frac{1}{p}} = |f|_{B_{p,\infty}^\alpha} < +\infty,$$

and it belongs to $B_{p,q}^\alpha([0,1]^l)$ with $q < +\infty$ if

$$\sum_{j \geq 0} \left[ 2^{j(\alpha + \frac{l}{2} - \frac{l}{p})} \left( \sum_{\boldsymbol{k} \in \Lambda(j)} |\beta_{j,\boldsymbol{k}}|^p \right)^{\frac{1}{p}} \right]^q = |f|_{B_{p,q}^\alpha}^q < +\infty.$$

Many properties of those function spaces are to be found in DeVore and Lorentz [19], DeVore [17] and Härdle, Kerkyacharian, Picard and Tsybakov [22] among other references.

As a consequence of Theorem 1, we can derive an adaptation result for the estimation of the intensity of a Poisson process when it belongs to some Besov space on $[0,1]^l$.

**Theorem 4.** *Let $\boldsymbol{X}$ be a Poisson process with unknown intensity $s$ with respect to Lebesgue measure on $[0,1]^l$. Let us assume that $\sqrt{s}$ belongs to some Besov space $B_{p,\infty}^\alpha([0,1]^l)$ for some unknown values of $p > 0$, $\alpha > l(1/p - 1/2)_+$ and $|\sqrt{s}|_{B_{p,\infty}^\alpha}$ given by (4.1). One can build a T-estimator $\hat{s}(\boldsymbol{X})$ such that*

$$(4.2) \quad \mathbb{E}_s\left[H^2(s, \hat{s})\right] \leq C(\alpha, p, l) \left[|\sqrt{s}|_{B_{p,\infty}^\alpha} \vee 1\right]^{2l/(2\alpha + l)}.$$

*Proof.* We just use Proposition 13 of Birgé [9] which provides suitable families $\mathcal{M}_j(2^i)$ of linear approximation spaces for functions in $B_{p,\infty}^\alpha([0,1]^l)$ and use the family of linear spaces $\{\overline{S}_m\}_{m \in \mathcal{M}}$ with $\mathcal{M} = \bigcup_{i \geq 1} \bigcup_{j \geq 0} \mathcal{M}_j(2^i)$ provided by this proposition. Then, for $m \in \mathcal{M}_j(2^i)$, $\overline{D}_m \leq c_1(2^i) + c_2(2^i)2^{jl}$ and we choose $\Delta_m = c_3(2^i)2^{jl} + i + j$ which implies that (1.13) holds with $\Sigma < 1$. Applying Proposition 13 of Birgé [9] with $t = \sqrt{s}$, $r = 2^i > \alpha \geq 2^{i-1}$ and $q = 2$, we derive from Theorem 1 that, if $R = |\sqrt{s}|_{B_{p,\infty}^\alpha} \vee 1$,

$$\mathbb{E}_s\left[H^2(s, \hat{s})\right] \leq C \inf_{j \geq 0}\left\{C(\alpha, p, l)R^2 2^{-2j\alpha} + c_4(\alpha)2^{jl}\right\}.$$

Choosing for $j$ the smallest integer such that $2^{j(l+2\alpha)} \geq R^2$ leads to the result. □

### *4.2. Anisotropic Hölder spaces*

Let us recall that a function $f$ defined on $[0,1)$ belongs to the Hölder class $\mathcal{H}(\alpha, R)$ with $\alpha = \beta + p$, $p \in \mathbb{N}$, $0 < \beta \leq 1$ and $R > 0$ if $f$ has a derivative of order $p$ satisfying $|f^{(p)}(x) - f^{(p)}(y)| \leq R|x-y|^\beta$ for all $x, y \in [0, 1)$. Given two multi-indices $\boldsymbol{\alpha} = (\alpha_1, \ldots, \alpha_k)$ and $\boldsymbol{R} = (R_1, \ldots, R_k)$ in $(0, +\infty)^k$, we define the anisotropic Hölder class $\mathcal{H}(\boldsymbol{\alpha}, \boldsymbol{R})$ as the set of functions $f$ on $[0,1)^k$ such that, for each $j$ and each set of $k-1$ coordinates $x_1, \ldots, x_{j-1}, x_{j+1}, \ldots, x_k$ the univariate function $y \mapsto f(x_1, \ldots, x_{j-1}, y, x_{j+1}, \ldots, x_k)$ belongs to $\mathcal{H}(\alpha_j, R_j)$.

Let now a multi-integer $\boldsymbol{N} = (N_1, \ldots, N_k) \in (\mathbb{N}^\star)^k$ be given. To it corresponds the hyperrectangle $\prod_{j=1}^k [0, N_j^{-1})$ and the partition $\mathcal{I}_{\boldsymbol{N}}$ of $[0,1)^k$ into $\prod_{j=1}^k N_j$ translates of this hyperrectangle. Given an integer $r \in \mathbb{N}$ and $m = (\boldsymbol{N}, r)$ we can define the linear space $\overline{S}_m$ of piecewise polynomials on the partition $\mathcal{I}_{\boldsymbol{N}}$ with degree at most $r$ with respect to each variable. Its dimension is $\overline{D}_m = (r+1)^k \prod_{j=1}^k N_j$. Setting $\mathcal{M} = (\mathbb{N}^\star)^k \times \mathbb{N}$ and $\Delta_m = \overline{D}_m$, we get (1.13) with $\Sigma$ depending only on $k$ as shown in the proof of Proposition 5, page 346 of Barron, Birgé and Massart [7]. The same proof also implies (see (4.25), page 347) the following approximation lemma.

**Lemma 4.** *Let $f \in \mathcal{H}(\boldsymbol{\alpha}, \boldsymbol{R})$ with $\alpha_j = \beta_j + p_j$, $r \geq \max_{1 \leq j \leq k} p_j$, $\boldsymbol{N} = (N_1, \ldots, N_k) \in (\mathbb{N}^\star)^k$ and $m = (\boldsymbol{N}, r)$. There exists some $g \in \overline{S}_m$ such that*

$$\|f - g\|_\infty \leq C(k, r) \sum_{j=1}^k R_j N_j^{-\alpha_j}.$$

We are now in a position to state the following corollary of Theorem 1.

**Corollary 2.** *Let $\boldsymbol{X}$ be a Poisson process with unknown intensity $s$ with respect to the Lebesgue measure on $[0,1)^k$ and $\hat{s}$ be a T-estimator based on the family of linear models $\{\overline{S}_m, m \in \mathcal{M}\}$ that we have previously defined. Assume that $\sqrt{s}$ belongs to the class $\mathcal{H}(\boldsymbol{\alpha}, \boldsymbol{R})$ and set*

$$\overline{\alpha} = \left(k^{-1} \sum_{j=1}^k \alpha_j^{-1}\right)^{-1} \quad \text{and} \quad \overline{R} = \left(\prod_{j=1}^k R_j^{1/\alpha_j}\right)^{\overline{\alpha}/k}.$$

*If $R_j \geq \overline{R}^{k/(2\overline{\alpha}+k)}$ for all $j$, then*

$$\mathbb{E}_s\left[H^2(s, \hat{s})\right] \leq C(k, \boldsymbol{\alpha})\overline{R}^{2k/(2\overline{\alpha}+k)}.$$

*Proof.* If $\alpha_j = \beta_j + p_j$ for $1 \leq j \leq k$, let us set $r = \max_{1 \leq j \leq k} p_j$, $\eta = \overline{R}^{k/(2\overline{\alpha}+k)}$ and define $N_j \in \mathbb{N}^\star$ by $(R_j/\eta)^{1/\alpha_j} \leq N_j < (R_j/\eta)^{1/\alpha_j} + 1$ so that $N_j < 2(R_j/\eta)^{1/\alpha_j}$ for all $j$. It follows from Lemma 4 that there exists some $t \in \overline{S}_m$, $m = (\boldsymbol{N}, r)$ with $\|\sqrt{s} - t\|_\infty \leq C_1(k, \boldsymbol{\alpha}) \sum_{j=1}^k R_j N_j^{-\alpha_j}$, hence $\|\sqrt{s} - t\|_2 \leq kC_1(k, \boldsymbol{\alpha})\eta$. It then follows from Theorem 1 that

$$\mathbb{E}_s\left[H^2(s, \hat{s})\right] \leq C_2(k, \boldsymbol{\alpha}) \left[\eta^2 + (r+1)^k \prod_{j=1}^k N_j\right] \leq C_3(k, \boldsymbol{\alpha}) \left[\eta^2 + \overline{R}^{k/\overline{\alpha}} \eta^{-k/\overline{\alpha}}\right].$$

The conclusion follows. □

### 4.3. Intensities with bounded $\alpha$-variation

Let us first recall that a function $f$ defined on some interval $J \subset \mathbb{R}$ has bounded $\alpha$-variation on $J$ for some $\alpha \in (0, 1]$ if

$$(4.3) \qquad \sup_{i \geq 1} \sup_{\substack{x_0 < \cdots < x_i \\ x_j \in J \text{ for } 0 \leq j \leq i}} \sum_{j=1}^i |f(x_j) - f(x_{j-1})|^{1/\alpha} = [V_\alpha(f; J)]^{1/\alpha} < +\infty,$$

the classical case of bounded variation corresponding to $\alpha = 1$. This formulation using the power $1/\alpha$ (instead of $\alpha$) implies that an $\alpha$-Hölderian function has bounded $\alpha$-variation over any finite interval $J$. We want to build a family of linear models which are suitable for estimating intensities $s$ with support on some interval $J$ of finite length $L$ and such that $\sqrt{s}$ has bounded $\alpha$-variation on $J$ for some unknown value of $\alpha$. These models are linear spaces of piecewise constant functions on some finite partitions $m$ of $J$, namely

$$\overline{S}_m = \left\{t = \sum_{j=1}^D a_j \mathbb{1}_{I_j}\right\} \quad \text{when } m = \{I_1, \ldots, I_D\}.$$

We consider for $\mathcal{M}$ a special family of partitions $m$ of $J$ derived by dyadic splitting which are in one-to-one correspondence with the family of complete binary trees. They are built according to the following "adaptive" algorithm described in Section 3.3 of DeVore [17]. This algorithm simultaneously grows a complete binary tree and a dyadic partition of $J$. It starts with a tree reduced to its root which is associated to the interval $J$. At each step of the algorithm the set of terminal nodes of the current tree is associated to the set of intervals in the current partition. Each step of the algorithm corresponds to choosing one terminal node and adding two sons to it. For the associated partition this means dividing the interval which corresponds to this terminal node into two intervals of equal length which then correspond to the two sons. At some stage the procedure stops and we end with a complete binary tree with $D$ terminal nodes and the associated partition of $J$ into $D$ intervals. We acually take for $\mathcal{M}$ the set of all finite partitions $m$ that can be build in that way so that each $m$ corresponds to the complete binary tree with $|m|$ terminal nodes that was used to build the partition.

It is known that the number of complete binary trees with $j + 1$ terminal nodes is given by the so-called Catalan numbers $(1+j)^{-1}\binom{2j}{j} \leq 4^j/(1+j)$ as explained for instance in Stanley [33], page 172. Setting $\Delta_m = 2|m|$ leads to

$$\sum_{m \in \mathcal{M}} \exp[-\Delta_m] = \sum_{j \geq 0} \sum_{\{m \in \mathcal{M} \mid |m|=1+j\}} \exp[-2(j+1)]$$

(4.4)
$$\leq \sum_{j \geq 0} \frac{4^j \exp[-2(j+1)]}{j+1} = e^{-2} \sum_{j \geq 0} \frac{(2/e)^{2j}}{j+1} < 1.$$

The approximation properties of $\bigcup_{m \in \mathcal{M}} \overline{S}_m$ with respect to functions of bounded $\alpha$-variation are given by the following proposition the proof of which was kindly communicated to the author by Ron DeVore [18].

**Proposition 3.** *Let $f$ be a function of bounded $\alpha$-variation on the interval $J$ of finite length $L$ with $\alpha$-variation $V_\alpha(f;J)$ given by (4.3). For each $j \in \mathbb{N}$, one can find a partition $m \in \mathcal{M}$ with*

(4.5) $\qquad |m| \leq c_1(\alpha) 2^j \qquad and \qquad \inf_{t \in \overline{S}_m} \|f-t\|_2 \leq c_2(\alpha) L^{1/2} V_\alpha(f;J) 2^{-j\alpha}.$

*with $1 < c_1(\alpha) = (1 - 2^{-[1/(2\alpha)+1]})(1 - 2^{-1/(2\alpha)})^{-1} < 2.21$ and*

$$\sqrt{2} < c_2(\alpha) = \left[\frac{2^{1+2\alpha}\left(1 - 2^{-[1/(2\alpha)+1]}\right)^{1-2\alpha}}{1 - 2^{-1/(2\alpha)}}\right]^{1/2} < 6.51.$$

*Proof.* For any interval $I \subset J$ we denote by $|I|$ its length and set $V(I) = V_\alpha(f;I)$. If $m = \{I_1; \ldots; I_D\}$ is a partition of $J$ into $D$ intervals, $\bar{f}_j = |I_j|^{-1} \int_{I_j} f(x)\,dx$ and $\bar{f} = \sum_{j=1}^D \bar{f}_j \mathbb{1}_{I_j}$, then $\|(f - \bar{f}_j)\mathbb{1}_{I_j}\|_\infty \leq V(I_j)$, hence

(4.6) $\qquad \|f - \bar{f}\|_2^2 \leq \sum_{j=1}^D E(I_j) \quad \text{with} \quad E(I) = |I| V^2(I).$

In particular (4.5) holds with $m = \{J\}$ and $j = 0$. To study the general case we choose some $\varepsilon > 0$ and apply the adaptive algorithm described just before in the following way: at each step we inspect the intervals of the partition and if we find an interval $I$ with $E(I) > \varepsilon$ we divide it into two intervals of equal length $|I|/2$. The algorithm necessarily stops since $E(I) \leq |I| V^2(J)$ for all $I \subset J$ and this results in some partition $m$ with $E(I) \leq \varepsilon$ for all $I \in m$. It follows from (4.6) that if $\bar{f}$ is built on this partition, then $\|f - \bar{f}\|_2^2 \leq \varepsilon |m|$. Since the case $|m| = 1$ has already been considered, we may assume that $|m| \geq 2$. Let us denote by $D_k$ the number of intervals in $m$ with length $L2^{-k}$ and set $a_k = 2^{-k} D_k$ so that $\sum_{k \geq 1} a_k = 1$ (since $D_0 = 0$). If $I$ is an interval of length $L2^{-k}$, $k > 0$, it derives from the splitting of an interval $I'$ with length $L2^{-k+1}$ such that $E(I') > \varepsilon$, hence, by (4.6), $V(I') > [\varepsilon L^{-1} 2^{k-1}]^{1/2}$ and, since the set function $V^{1/\alpha}$ is subadditive over disjoint intervals, the number of such interval $I'$ is bounded by $[V(J)]^{1/\alpha} [\varepsilon L^{-1} 2^{k-1}]^{-1/(2\alpha)}$. It follows that

$$D_k \leq \gamma 2^{-k/(2\alpha)} \quad \text{and} \quad a_k \leq \gamma 2^{-k/(2\alpha)-k} \quad \text{with} \quad \gamma = 2[V(J)]^{1/\alpha} [\varepsilon/(2L)]^{-1/(2\alpha)}.$$

Since $|m| = \sum_{k \geq 1} 2^k a_k$, we can derive a bound on $|m|$ from a maximization of $\sum_{k \geq 1} 2^k a_k$ under the restrictions $\sum_{k \geq 1} a_k = 1$ and $a_k \leq \gamma 2^{-k[1/(2\alpha)+1]}$. One should then clearly keep the largest possible indices $k$ with the largest possible values for $a_k$. Let us fix $\varepsilon$ so that $\gamma = (1 - 2^{-[1/(2\alpha)+1]})2^{j[1/(2\alpha)+1]}$ for some $j \geq 1$. Then, setting $a_k$ to its maximal value, we get $\sum_{k \geq j} \gamma 2^{-k[1/(2\alpha)+1]} = 1$, which implies that an upper bound for $|m|$ is

$$|m| \leq \sum_{k \geq j} \gamma 2^k 2^{-k[1/(2\alpha)+1]} = \frac{\gamma 2^{-j/(2\alpha)}}{1 - 2^{-1/(2\alpha)}} = \frac{1 - 2^{-[1/(2\alpha)+1]}}{1 - 2^{-1/(2\alpha)}} 2^j.$$

The corresponding value of $\varepsilon$ is $2L(\gamma/2)^{-2\alpha}V^2(J)$ so that

$$\|f - \bar{f}\|_2^2 \leq \varepsilon |m| \leq 2LV^2(J)2^{2\alpha} \frac{\gamma^{1-2\alpha} 2^{-j/(2\alpha)}}{1 - 2^{-1/(2\alpha)}}$$
$$= \frac{2LV^2(J)2^{2\alpha} \left(1 - 2^{-[1/(2\alpha)+1]}\right)^{1-2\alpha}}{1 - 2^{-1/(2\alpha)}} 2^{-2\alpha j}.$$

These two bounds give (4.5) and we finally use the fact that $0 < \alpha \leq 1$ to bound the two constants. $\square$

We can then derive from this proposition, (1.15) and our choice of the $\Delta_m$ that

$$\mathbb{E}_s \left[ H^q(s, \hat{s}) \right] \leq C(q) \inf_{j \in \mathbb{N}} \left\{ 2^{j/2} + L^{1/2} V_\alpha \left( \sqrt{s}; J \right) 2^{-j\alpha} \right\}^q.$$

An optimization with respect to $j \in \mathbb{N}$ then leads to the following risk bound.

**Corollary 3.** *Let $\boldsymbol{X}$ be a Poisson process with unknown intensity $s$ with respect to the Lebesgue measure on some interval $J$ of length $L$. We assume that $\sqrt{s}$ has finite $\alpha$-variation equal to $V$ on $J$, both $\alpha$ and $V$ being unknown. One can build a T-estimator $\hat{s}(\boldsymbol{X})$ such that*

(4.7) $$\mathbb{E}_s \left[ H^q(s, \hat{s}) \right] \leq C(q) \left[ \left( L^{1/2} V \right) \vee 1 \right]^{q/(2\alpha+1)}.$$

It is not difficult to show, using Assouad's Lemma, that, up to a constant, this bound is optimal when $q = 2$.

**Proposition 4.** *Let $L, \alpha$ and $V$ be given and $\mathcal{S} \subset \mathbb{L}_1^+(\lambda)$ be the set of intensities with respect to the Lebesgue measure on $[0, L]$ such that $\sqrt{s}$ has $\alpha$-variation bounded by $V$. Let $\hat{s}(\boldsymbol{X})$ be any estimator based on a Poisson process $\boldsymbol{X}$ with unknown intensity $s \in \mathcal{S}$. There exists a universal constant $c > 0$ (independent of $\hat{s}, L, \alpha$ and $V$) such that*

$$\sup_{s \in \mathcal{S}} \mathbb{E}_s \left[ H^2(s, \hat{s}) \right] \geq c \left[ \left( L^{1/2} V \right) \vee 1 \right]^{2/(2\alpha+1)}.$$

*Proof.* If $L^{1/2}V < 1$, we simply apply (2.6) with $s_0 = \mathbb{1}_{[0,L)}$ and $s_1 = (1 + L^{-1/2})^2 \mathbb{1}_{[0,L)}$ so that $2H^2(s_0, s_1) = 1$. If $L = 1$ and $V \geq 1$ we fix some positive integer $D$ and define $g$ with support on $[0, D^{-1})$ by

$$g(x) = x \mathbb{1}_{[0,(2D)^{-1})}(x) + \left( D^{-1} - x \right) \mathbb{1}_{[(2D)^{-1}, D^{-1})}(x).$$

Then $\int_0^{1/D} g^2(x)\, dx = (12D^3)^{-1}$ and $0 \leq g(x) \leq (2D)^{-1}$. If we apply the construction of Lemma 2, we get a family of Lipschitz intensities $s_\delta$ with values in the interval $[12D^3 - 3D^2, 12D^3 + 3D^2] \subset [9D^3, 15D^3]$ and Lipschitz coefficient $6D^3$. It follows that if $0 \leq x < y \leq 1$,

$$\left| \sqrt{s_\delta(x)} - \sqrt{s_\delta(y)} \right| \leq \frac{|s_\delta(x) - s_\delta(y)|}{6D^{3/2}}$$
$$\leq \frac{(6D^2) \wedge (6D^3|x-y|)}{6D^{3/2}} \leq \sqrt{D}\left[1 \wedge (D|x-y|)\right].$$

This allows us to bound the $\alpha$-variation of $\sqrt{s_\delta}$ in the following way. For any increasing sequence $0 \leq x_0 < \cdots < x_i \leq 1$,

$$\sum_{j=1}^{i} \left|\sqrt{s_\delta(x_j)} - \sqrt{s_\delta(x_{j-1})}\right|^{1/\alpha} \leq D^{1/(2\alpha)} \sum_{j=1}^{i} \mathbb{1}_{\{x_j - x_{j-1} \geq D^{-1}\}}$$
$$+ D^{3/(2\alpha)} \sum_{j=1}^{i} \mathbb{1}_{\{x_j - x_{j-1} < D^{-1}\}} (x_j - x_{j-1})^{1/\alpha}.$$

If $n = \sum_{j=1}^{i} \mathbb{1}_{\{x_j - x_{j-1} \geq D^{-1}\}} \leq D$, then

$$D^{3/(2\alpha)} \sum_{j=1}^{i} \mathbb{1}_{\{x_j - x_{j-1} < D^{-1}\}} (x_j - x_{j-1})^{1/\alpha}$$
$$\leq D^{3/(2\alpha)} D^{-1/\alpha} (D - n) = D^{1/(2\alpha)}(D - n),$$

which shows that the $\alpha$-variation of $\sqrt{s_\delta}$ is bounded by $[D^{1/(2\alpha)}D]^\alpha = D^{(1+2\alpha)/2}$. We finally choose for $D$ the largest integer $j$ such that $j^{(1+2\alpha)/2} \leq V$. Then $V^{2/(1+2\alpha)} < 2D$ and an application of Lemmas 1 and 2 show that

$$\sup_{s \in \mathcal{S}_D} \mathbb{E}_s \left[ H^2(s, \hat{s}) \right] \geq 2^{-8}(2D) \exp[-2/7] \geq 2^{-8} \exp[-2/7] V^{2/(1+2\alpha)},$$

which proves our lower bound. The general case $L^{1/2}V \geq 1$ follows from a scaling argument. If $\boldsymbol{X}$ is a Poisson process on $[0, L]$ with intensity $s$ (with respect to the Lebesgue measure), then $\boldsymbol{Y} = L^{-1}\boldsymbol{X}$ is a Poisson process on $[0, 1]$ with intensity $s_L$ to which the previous results apply. Since $s_L(y) = Ls(Ly)$, it follows that $H^2(s, t) = H^2(s_L, t_L)$ and, if $\sqrt{s}$ has $\alpha$-variation bounded by $V$, $\sqrt{s_L}$ has $\alpha$-variation bounded by $L^{1/2}V$. The result for an arbitrary $L$ follows from these remarks. $\square$

### 4.4. Intensities with square roots in weak $\ell_q$-spaces

#### 4.4.1. Approximation based on weak $\ell_q$-spaces

As we already mentioned, if $s \in \mathbb{L}_1^+(\lambda)$ is an intensity with respect to $\lambda$ on $\mathcal{X}$ and we are given an orthonormal basis $\{\varphi_j, j \geq 1\}$ of $\mathbb{L}_2(\lambda)$, $\sqrt{s}$ can be written as $\sum_{j \geq 1} \beta_j \varphi_j$ with $\boldsymbol{\beta} = (\beta_j)_{j \geq 1} \in \ell_2 = \ell_2(\mathbb{N}^*)$ and $\sum_{j \geq 1} \beta_j^2 = \|\sqrt{s}\|_2^2 < +\infty$. Hence, for all $x > 0$, $|\{j \geq 1 \,|\, |\beta_j| \geq x\}| \leq \|\sqrt{s}\|_2^2 x^{-2}$, which means that the sequence $\boldsymbol{\beta}$ belongs to the weak $\ell_2$-space $\ell_2^w$.

More generally, given a sequence $\boldsymbol{\beta} = (\beta_j)_{j \geq 1}$ converging to zero and $a_j$ the rearrangement of the numbers $|\beta_j|$ in nonincreasing order (which means that $a_1 = \sup_{j \geq 1} |\beta_j|$, etc...), we say that $\boldsymbol{\beta}$ belongs to the weak $\ell_q$-space $\ell_q^w$ $(q > 0)$ if

(4.8) $\quad \sup_{x > 0} x^q \,|\{j \geq 1 \,|\, |\beta_j| \geq x\}| = \sup_{x > 0} x^q \,|\{j \geq 1 \,|\, a_j \geq x\}| = |\boldsymbol{\beta}|_{q,w}^q < +\infty.$

This implies that $a_j \leq |\boldsymbol{\beta}|_{q,w} j^{-1/q}$ for $j \geq 1$ and the reciprocal actually holds:

(4.9) $\quad\quad |\boldsymbol{\beta}|_{q,w} = \inf \left\{ y > 0 \,|\, a_j \leq y j^{-1/q} \quad \text{for all } j \geq 1 \right\}.$

Note that, although $|\theta\boldsymbol{\beta}|_{q,w} = |\theta||\boldsymbol{\beta}|_{q,w}$ for $\theta \in \mathbb{R}$, $|\boldsymbol{\beta}|_{q,w}$ is not a norm. For convenience, we shall call it the *weight* of $\boldsymbol{\beta}$ in $\ell_q^w$. By extension, given the basis

$\{\varphi_j, j \geq 1\}$, we shall say that $u \in \mathbb{L}_2(\lambda)$ belongs to $\ell_q^w$ if $u = \sum_{j \geq 1} \beta_j \varphi_j$ and $\boldsymbol{\beta} \in \ell_q^w$. As a consequence of this control on the size of the coefficients $a_j$, we get the following useful lemma.

**Lemma 5.** *Let $\boldsymbol{\beta} \in \ell_q^w$ with weight $|\boldsymbol{\beta}|_{q,w}$ for some $q > 0$ and $(a_j)_{j \geq 1}$ be the nonincreasing rearrangement of the numbers $|\beta_j|$. Then $\boldsymbol{\beta} \in \ell_p$ for $p > q$ and for all $n \geq 1$,*

$$\text{(4.10)} \qquad \sum_{j > n} a_j^p \leq \frac{q}{p - q} |\boldsymbol{\beta}|_{q,w}^p (n + 1/2)^{-(p-q)/q}.$$

*Proof.* By (4.9) and convexity,

$$\sum_{j>n} a_j^p \leq |\boldsymbol{\beta}|_{q,w}^p \sum_{j>n} j^{-p/q} \leq |\boldsymbol{\beta}|_{q,w}^p \int_{n+1/2}^{+\infty} x^{-p/q}\, dx. \qquad \square$$

As explained in great detail in Kerkyacharian and Picard [23] and Cohen, DeVore, Kerkyacharian and Picard [15], the fact that $u \in \ell_q^w$ for some $q < 2$ has important consequences for the approximation of $u$ by fonctions in suitable $D$-dimensional spaces. For $m$ any finite subset of $\mathbb{N}^\star$, let us define $\overline{S}_m$ as the linear span of $\{\varphi_j, j \in m\}$. If $u = \sum_{j \geq 1} \beta_j \varphi_j$ belongs to $\ell_q^w$ and $D$ is a positive integer, one can find some $m$ with $|m| = D$ and some $t \in \overline{S}_m$ such that

$$\text{(4.11)} \qquad \|u - t\|_2^2 \leq (2/q - 1)^{-1} |\boldsymbol{\beta}|_{q,w}^2 (D + 1/2)^{1 - 2/q}.$$

Indeed, let us take for $m$ the set of indices of the $D$ largest numbers $|\beta_j|$. It follows from (4.10) that

$$\sum_{j \notin m} \beta_j^2 = \sum_{j > D} a_j^2 \leq \frac{q}{2 - q} |\boldsymbol{\beta}|_{q,w}^2 (D + 1/2)^{1-2/q}.$$

Setting $t = \sum_{j \in m} \beta_j \varphi_j$ gives (4.11) which provides the rate of approximation of $u$ by functions of the set $\bigcup_{\{m \,|\, |m|=D\}} \overline{S}_m$ as a decreasing function of $D$ (which is not possible for $q = 2$). Unfortunately, this involves an infinite family of linear spaces $\overline{S}_m$ of dimension $D$ since the largest coefficients of the sequence $\boldsymbol{\beta}$ may have arbitrarily large indices. To derive a useful, as well as a practical approximation method for functions in $\ell_q^w$-spaces, one has to restrict to those sets $m$ which are subsets of $\{1, \ldots, n\}$ for some given value of $n$. This is what is done in Kerkyacharian and Picard [23] who show, in their Corollary 3.1, that a suitable thresholding of empirical versions of the coefficients $\beta_j$ for $j \in \{1, \ldots, n\}$ leads to estimators that have nice properties. Of course, since this approach ignores the (possibly large) coefficients with indices bigger than $n$, an additional condition on $\boldsymbol{\beta}$ is required to control $\sum_{j > n} \beta_j^2$. In Kerkyacharian and Picard [23], it takes the form

$$\text{(4.12)} \qquad \sum_{j > n} \beta_j^2 \leq A^2 n^{-\delta} \quad \text{for all } n \geq 1, \quad \text{with } A \text{ and } \delta > 0,$$

while Cohen, DeVore, Kerkyacharian and Picard [15], page 178, use the similar condition BS. Such a condition is always satisfied for functions in Besov spaces $B_{p,\infty}^\alpha([0,1]^l)$ with $p \leq 2$ and $\alpha > l(1/p - 1/2)$. Indeed, if

$$f = \sum_{j=-1}^\infty \sum_{\boldsymbol{k} \in \Lambda(j)} \beta_{j,\boldsymbol{k}} \varphi_{j,\boldsymbol{k}}$$

belongs to such a Besov space, it follows from (4.1) that,

$$\sum_{j>J}\sum_{\mathbf{k}\in\Lambda(j)}|\beta_{j,\mathbf{k}}|^2 \leq \sum_{j>J}\left(\sum_{\mathbf{k}\in\Lambda(j)}|\beta_{j,\mathbf{k}}|^p\right)^{2/p} \leq |f|^2_{B^\alpha_{p,\infty}}\sum_{j>J}2^{-2j\left(\alpha+\frac{1}{2}-\frac{l}{p}\right)}$$

(4.13)
$$\leq C|f|^2_{B^\alpha_{p,\infty}}2^{-2J\left(\alpha+\frac{1}{2}-\frac{l}{p}\right)}.$$

Since the number of coefficients $\beta_{j,\mathbf{k}}$ with $j \leq J$ is bounded by $C'2^{Jl}$, after a proper change in the indexing of the coefficients, the corresponding sequence $\boldsymbol{\beta}$ will satisfy $\sum_{j>n}\beta_j^2 \leq A^2 n^{-\delta}$ with $\delta = (2\alpha/l) + 1 - (2/p)$.

### 4.4.2. Model selection for weak $\ell_q$-spaces

It is the very method of thresholding that imposes to fix the value of $n$ as a function of $\delta$ or impose the value of $\delta$ when $n$ has been chosen in order to get a good performance for the threshold estimators. Model selection is more flexible since it allows to adapt the value of $n$ to the unknown values of $A$ and $\delta$. Let us assume that an orthonormal basis $\{\varphi_j, j \geq 1\}$ for $\mathbb{L}_2(\lambda)$ has been chosen and that the Poisson process $X$ has an intensity $s$ with respect to $\lambda$ so that $\sqrt{s} = \sum_{j\geq 1}\beta_j\varphi_j$ with $\boldsymbol{\beta} \in \ell_2$. We take for $\mathcal{M}$ the set of all subsets $m$ of $\mathbb{N}^\star$ such that $|m| = 2^j$ for some $j \in \mathbb{N}$ and choose for $\overline{S}_m$ the linear span of $\{\varphi_j, j \in m\}$ with dimension $\overline{D}_m = |m|$. If $|m| = 2^j$ and $k = \inf\{i \in \mathbb{N}^\star \,|\, 2^i \geq l \quad \text{for all } l \in m\}$, we set $\Delta_m = k + \log\binom{2^k}{2^j}$. Then

$$\sum_{m\in\mathcal{M}}\exp[-\Delta_m] \leq \sum_{k\geq 1}\sum_{j=0}^k\binom{2^k}{2^j}\exp\left[-k - \log\binom{2^k}{2^j}\right] \leq \sum_{k\geq 1}(k+1)\exp[-k],$$

which allows to apply Theorem 1.

**Proposition 5.** *Let $\hat{s}$ be a T-estimator provided by Theorem 1 and based on the previous family of models $\overline{S}_m$ and weights $\Delta_m$. If $\sqrt{s} = \sum_{j\geq 1}\beta_j\varphi_j$ with $\boldsymbol{\beta} \in \ell_q^w$ for some $q < 2$ and (4.12) holds with $A \geq 1$ and $0 < \delta \leq 1$, the risk of $\hat{s}$ at $s$ is bounded by*

$$\mathbb{E}_s\left[H^2(s,\hat{s})\right] \leq C\left[\left(\gamma^{1-q/2}\left(R^2 \vee \gamma\right)^{q/2}\right) \bigwedge A^{2/(1+\delta)}\right],$$

with
$$R = \left[\frac{q}{2-q}\right]^{1/2}|\boldsymbol{\beta}|_{q,w} \quad \text{and} \quad \gamma = \delta^{-1}\left[\frac{\log\left(\delta[A \vee R]^2\right)}{\log 2}\bigvee 1\right].$$

*Proof.* Let $(a_j)_{j\geq 1}$ be the nonincreasing rearrangement of the numbers $|\beta_j|$, $k$ and $j \leq k$ be given and $m$ be the set of indices of the $2^j$ largest coefficients among $\{|\beta_1|,\ldots,|\beta_{2^k}|\}$. Then $\overline{D}_m = 2^j$ and $\Delta_m \leq k + \log\binom{2^k}{2^j}$. It follows from (4.10) and (4.12) that

$$\sum_{j\notin m}\beta_j^2 \leq \left(\sum_{i>2^j}a_i^2\right)\mathbb{1}_{j<k} + \sum_{i>2^k}\beta_i^2 \leq \frac{q}{2-q}|\boldsymbol{\beta}|^2_{q,w}2^{-j(2/q-1)}\mathbb{1}_{j<k} + A^2 2^{-k\delta}.$$

This shows that one can find $t \in \overline{S}_m$ such that $\|\sqrt{s} - t\|_2^2 \leq R^2 2^{-j(2/q-1)}\mathbb{1}_{j<k} + A^2 2^{-k\delta}$ and it follows from (1.14) that

$$\mathbb{E}_s\left[H^2(s,\hat{s})\right] \leq C\inf_{k\geq 1}\inf_{0\leq j\leq k}\left\{R^2 2^{-j(2/q-1)}\mathbb{1}_{j<k} + A^2 2^{-k\delta} + 2^j + k + \log\binom{2^k}{2^j}\right\}.$$

We recall that $C$ denotes a constant that may change as often as necessary. If $j = k$, $\mathbb{E}_s[H^2(s,\hat{s})] \leq C[A^2 2^{-k\delta} + 2^k]$ and an optimization with respect to $k$ leads to $\mathbb{E}_s[H^2(s,\hat{s})] \leq CA^{2/(1+\delta)}$. For $j < k$, we notice that $\Delta_m \leq k + 2^j[1+\log(2^{k-j})] < 3k2^j$, so that

$$(4.14) \quad \mathbb{E}_s\left[H^2(s,\hat{s})\right] \leq C \inf_{k \geq 1}\left\{\left(A^2 2^{-k\delta}\right) \vee \inf_{0 \leq j < k}\left\{\left(R^2 2^{-j(2/q-1)}\right) \vee (k2^j)\right\}\right\}.$$

If $R^2 2^{-(k-1)(2/q-1)} > k2^{k-1}$, we may harmlessly increase $k$ until $k = K$ with

$$K = \inf\left\{i \geq 1 \,\Big|\, i2^{i-1} \geq R^2 2^{-(i-1)(2/q-1)}\right\} = \inf\left\{i \geq 1 \,\Big|\, 2^{i-1} \geq R^q i^{-q/2}\right\}$$

and therefore restrict the minimization in (4.14) to $k \geq K$. We then choose for $j$ the smallest integer $i$ such that $2^i \geq (R^2/k)^{q/2}$, which leads to

$$\mathbb{E}_s\left[H^2(s,\hat{s})\right] \leq C \inf_{k \geq K}\left\{\left(A^2 2^{-k\delta}\right) \vee \left(R^q k^{1-q/2}\right) \vee k\right\}.$$

It follows from Lemma 6 below (with $a = 1$) that, if $\delta A^2 \leq 2$, $(A^2 2^{-k\delta}) \bigvee k \geq A^2/2$ for all $k$ which does not improve on our previous bound $CA^{2/(1+\delta)}$ so that we may assume from now on that $\delta A^2 > 2$, hence $\gamma > \delta^{-1}$. Handling this case in full generality is much more delicate and we shall simplify the minimization problem by replacing $A$ by $\overline{A} = A \vee R$, which amounts to assuming that $A \geq R$ and leads to $\mathbb{E}_s[H^2(s,\hat{s})] \leq C \inf_{k \geq K} f(k)$ with

$$f(x) = f_1(x) \vee f_2(x) \vee x; \qquad f_1(x) = \overline{A}^2 2^{-x\delta} \quad \text{and} \quad f_2(x) = R^q x^{1-q/2}.$$

We want to minimize $f(x)$, up to constants. The minimization of $f_1(x) \vee x$ follows from Lemma 6 with $\delta \overline{A}^2 > 2$. The minimum then takes the form $c_2\gamma > 0.469\gamma$ with $f_1(\gamma) = \delta^{-1} < \gamma$ hence $f(\gamma) = \gamma \vee f_2(\gamma)$. To show that $\inf_x f(x) \geq cf(\gamma)$ when $\delta \overline{A}^2 > 2$, we distinguish between two cases. If $R^2 \leq \gamma$, $f(\gamma) = \gamma$ and we conclude from the fact that $\inf_x f(x) > 0.469\gamma$. If $R^2 > \gamma$, $f_2(x) > x$ for $x \leq \gamma$, $f(\gamma) = f_2(\gamma) > \gamma$ and the minimum of $f(x)$ is obtained for some $x_0 < \gamma$. Hence

$$\inf_x f(x) = \inf_x \{f_1(x) \vee f_2(x)\} = R^q \inf_x \left\{\left(B2^{-\delta x}\right) \vee x^{1-q/2}\right\} \quad \text{with } B = \overline{A}^2 R^{-q}.$$

It follows from Lemma 6 with $a = (2-q)/2$ that the result of this minimization depends on the value of

$$V = \frac{2\delta}{2-q}\overline{A}^{4/(2-q)} R^{-2q/(2-q)} = \frac{2\overline{A}^2 \delta}{2-q}\left(\frac{\overline{A}}{R}\right)^{2q/(2-q)} \geq \overline{A}^2 \delta > 2,$$

since $\overline{A} \geq R$. Then,

$$\inf_x f(x) \geq R^q \left[\frac{(2-q)\log V}{3\delta}\right]^{1-q/2} \geq R^q \gamma^{1-q/2}\left[\frac{(2-q)\log 2}{3}\right]^{1-q/2}$$
$$> 0.45 R^q \gamma^{1-q/2},$$

and we can conclude that, in both cases, $\inf_x f(x) \geq 0.45 f(\gamma)$. Let us now fix $k$ such that $\gamma + 1 \leq k < \gamma + 2$ so that $k < 3\gamma$. Then $2^{k-1} \geq 2^\gamma = (\overline{A}^2\delta)^{1/\delta}$ while $R^q k^{-q/2} \leq (R^2/\gamma)^{q/2} \leq (R^2\delta)^{q/2}$. This implies that $k \geq K$. Moreover $f(k) = k \vee f_2(k) < 3f(\gamma)$ which shows that $\inf_{k \geq K} f(k) < 3f(\gamma) < 6.7 \inf_x f(x)$ and justifies this choice of $k$. Finally $\mathbb{E}_s[H^2(s,\hat{s})] \leq C[\gamma \vee f_2(\gamma)]$. □

Note that our main assumption, namely that $\boldsymbol{\beta} \in \ell_q^w$, implies that $\sum_{j>n} a_j^p \leq R^2 n^{-2/q+1}$ by (4.10) while (4.12) entails that $\sum_{j>n} a_j^p \leq \sum_{j>n} \beta_j^p \leq A^2 n^{-\delta}$. Since it is only an additional assumption it should not be strictly stronger than the main one, which is the case if $A \leq R$ and $\delta \geq 2/q - 1$. It is therefore natural to assume that at least one of these inequalities does not hold.

**Lemma 6.** *For positive parameters $a, B$ and $\theta$, we consider on $\mathbb{R}_+$ the function $f(x) = B 2^{-\delta x} \vee x^a$. Let $V = a^{-1} \delta B^{1/a}$. If $V \leq 2$ then $\inf_x f(x) = c_1 B$ with $2^{-a} \leq c_1 < 1$. If $V > 2$, then $\inf_x f(x) = [c_2 a \delta^{-1} \log V]^a$ with $2/3 < c_2 < 1$.*

*Proof.* Clearly, the minimum is obtained when $x = x_0$ is the solution of $B 2^{-\delta x} = x^a$. Setting $x_0 = B^{1/a} y$ and taking base 2 logarithms leads to $y^{-1} \log_2(y^{-1}) = V$, hence $y < 1$. If $V \leq 2$, then $1 < y^{-1} \leq 2$ and the first result follows. If $V \geq 2$, the solution takes the form $y = z V^{-1} \log_2 V$ with $1 > z > [1 - (\log_2 V)^{-1} \log_2(\log_2 V)] > 0.469$. □

### 4.4.3. Intensities with bounded variation on $[0,1)^2$

This section, which is devoted to the estimation of an intensity $s$ such that $\sqrt{s}$ belongs to the space $BV([0,1)^2)$, owes a lot to discussions with Albert Cohen and Ron DeVore. The approximation results that we use here should be considered as theirs. The definition and properties of the space $BV([0,1)^2)$ of functions with bounded variation on $[0,1)^2$ are given in Cohen, DeVore, Petrushev and Xu [16] where the reader can also find the missing details. It is known that, with the notations of Section 4.1 for Besov spaces, $B_{1,1}^1([0,1)^2) \subset BV([0,1)^2) \subset B_{1,\infty}^1([0,1)^2)$. This corresponds to the situation $\alpha = 1, l = 2$ and $p = 1$, therefore $\alpha = l(1/p - 1/2)$, a borderline case which is not covered by the results of Theorem 4. On the other hand, it is proved in Cohen, DeVore, Petrushev and Xu [16], Section 8, that, if a function of $BV([0,1)^2)$ is expanded in the two-dimensional Haar basis, its coefficients belong to the space $\ell_1^w$. More precisely if $f \in BV([0,1)^2)$ with semi-norm $|f|_{BV}$ and $f$ is expanded in the Haar basis with coefficients $\beta_j$, then $|\boldsymbol{\beta}|_{1,w} \leq C |f|_{BV}$ where $|\boldsymbol{\beta}|_{1,w}$ is given by (4.8) and $C$ is a universal constant. We may therefore use the results of the previous section to estimate $\sqrt{s}$ but we need an additional assumption to ensure that (4.12) is satisfied. By definition $\sqrt{s}$ belongs to $\mathbb{L}_2([0,1)^2, dx)$ but we shall assume here slightly more, namely that it belongs to $\mathbb{L}_p([0,1)^2, dx)$ for some $p > 2$. This is enough to show that (4.12) holds.

**Lemma 7.** *If $f \in BV([0,1)^2) \cap \mathbb{L}_p([0,1)^2, dx)$ for some $p > 2$ and has an expansion $f = \sum_{j=-1}^{\infty} \sum_{\boldsymbol{k} \in \Lambda(j)} \beta_{j,\boldsymbol{k}} \varphi_{j,\boldsymbol{k}}$ with respect to the Haar basis on $[0,1)^2$, then for $J \geq -1$,*
$$\sum_{j>J} \sum_{\boldsymbol{k} \in \Lambda(j)} |\beta_{j,\boldsymbol{k}}|^2 \leq C(p) \|f\|_p |f|_{B_{1,\infty}^1} 2^{-2J(1/2 - 1/p)}.$$

*Proof.* It follows from Hölder inequality that $|\beta_{j,\boldsymbol{k}}| = \langle f, \varphi_{j,\boldsymbol{k}} \rangle \leq \|f\|_p \|\varphi_{j,\boldsymbol{k}}\|_{p'}$ with $p'^{-1} = 1 - p^{-1}$ and by the structure of a wavelet basis, $\|\varphi_{j,\boldsymbol{k}}\|_{p'}^{p'} \leq c_1 2^{-j(2-p')}$, so that $|\beta_{j,\boldsymbol{k}}| \leq c_2 \|f\|_p 2^{-j(2/p'-1)} = c_2 \|f\|_p 2^{-j(1-2/p)}$. Since $BV([0,1)^2) \subset B_{1,\infty}^1([0,1)^2)$, it follows from (4.1) with $\alpha = p = 1$ and $l = 2$ that $\sum_{\boldsymbol{k} \in \Lambda(j)} |\beta_{j,\boldsymbol{k}}| \leq |f|_{B_{1,\infty}^1}$ so that $\sum_{\boldsymbol{k} \in \Lambda(j)} |\beta_{j,\boldsymbol{k}}|^2 \leq c_2 \|f\|_p |f|_{B_{1,\infty}^1} 2^{-j(1-2/p)}$ for all $j \geq 0$. The conclusion follows. □

Since the number of coefficients $\beta_{j,\boldsymbol{k}}$ with $j \leq J$ is bounded by $C 2^{2J}$, after a proper reindexing of the coefficients, the corresponding sequence $\boldsymbol{\beta}$ will satisfy (4.12) with $\delta = 1/2 - 1/p$ which shows that it is essential here that $p$ be larger than 2. We finally get the following corollary of Proposition 5 with $q = 1$.

**Corollary 4.** *One can build a T-estimator $\hat{s}$ with the following properties. Let the intensity $s$ be such that $\sqrt{s} \in BV([0,1)^2) \cap \mathbb{L}_p([0,1)^2, dx)$ for some $p > 2$, so that the expansion of $\sqrt{s}$ in the Haar basis satisfies (4.12) with $\delta = 1/2 - 1/p$ and $A \geq 1$. Let $R = |\sqrt{s}|_{BV}$, then*

$$\mathbb{E}\left[H^2(s, \hat{s})\right] \leq C \left[\sqrt{\gamma(R^2 \vee \gamma)} \wedge A^{2/(1+\delta)}\right]$$

$$\text{with} \quad \gamma = \delta^{-1}\left[\frac{\log\left(\delta[A \vee R]^2\right)}{\log 2} \vee 1\right].$$

### 4.5. Mixing families of models

We have studied here a few families of approximating models. Many more can be considered and further examples can be found in Reynaud-Bouret [30] or previous papers of the author on model selection such as Barron, Birgé and Massart [7], Birgé and Massart [12], Birgé [9] and Baraud and Birgé [4]. As indicated in the previous sections, the choice of suitable families of models is driven by results in approximation theory relative to the type of intensity we expect to encounter or, more precisely, to the type of assumptions we make about the unknown function $\sqrt{s}$. Different types of assumptions will lead to different choices of approximating models, but it is always possible to combine them. If we have built a few families of linear models $\{\overline{S}_m, m \in \mathcal{M}_j\}$ for $1 \leq j \leq J$ and chosen suitable weights $\Delta_m$ such that $\sum_{m \in \mathcal{M}_j} \exp[-\Delta_m] \leq \Sigma$ for all $j$ we may consider the mixed family of models $\{\overline{S}_m, m \in \mathcal{M}\}$ with $\mathcal{M} = \cup_{j=1}^{J} \mathcal{M}_j$ and define new weights $\Delta'_m = \Delta_m + \log J$ for all $m \in \mathcal{M}$ so that (1.13) still holds with the same value of $\Sigma$. It follows from Theorem 1 that the T-estimator based on the mixed family will share the properties of the ones derived from the initial families apart, possibly, for a moderate increase in the risk of order $(\log J)^{q/2}$. The situation becomes more complex if $J$ is large or even infinite. A detailed discussion of how to mix families of models in general has been given in Birgé and Massart [12], Section 4.1, which applies with minor modifications to our case.

### 4.6. Asymptotics and a parallel with density estimation

The previous examples lead to somewhat unusual bounds with no number of observations $n$ like for density estimation and no variance size $\sigma^2$ as in the case of the estimation of a normal mean. Here, there is no rate of convergence because there is no sequence of experiments, just one with a mean measure $\mu_s = s \cdot \lambda$. To get back to more familiar results with rates and asymptotics and recover some classical risk bounds, we may reformulate our problem in a slightly different form which completely parallels the one we use for density estimation. As indicated in our introduction we may always rewrite the intensity $s$ as $s = ns_1$ with $\int s_1 \, d\lambda = 1$ so that $s_1$ becomes a density and $n = \mu_s(\mathcal{X})$. We use this notation here, although $n$ need not be an integer, to emphasize the similarity between the estimation of $s$ and density estimation. When $n$ is an integer this also corresponds to observing $n$ i.i.d. Poisson processes $\boldsymbol{X}_i$, $1 \leq i \leq n$ with intensity $s_1$ and set $\Lambda_{\boldsymbol{X}} = \sum_{i=1}^{n} \Lambda_{\boldsymbol{X}_i}$. In this case (1.15) can be rewritten in the following way.

**Corollary 5.** *Let $\lambda$ be some positive measure on $\mathcal{X}$, $\boldsymbol{X}$ be a Poisson process with unknown intensity $s \in \mathbb{L}_1^+(\lambda)$, $\{\overline{S}_m, m \in \mathcal{M}\}$ be a finite or countable family of linear subspaces of $\mathbb{L}_2(\lambda)$ with respective finite dimensions $\overline{D}_m$ and let $\{\Delta_m\}_{m \in \mathcal{M}}$ be a family of nonnegative weights satisfying (1.13). One can build a T-estimator $\hat{s}(\boldsymbol{X})$ of $s$ satisfying, for all $s \in \mathbb{L}_1^+(\lambda)$ such that $\int s\, d\lambda = n$, $s_1 = n^{-1}s$ and all $q \geq 1$,*

$$\mathbb{E}_s\left[\left(n^{-1/2}H(s,\hat{s})\right)^q\right] \leq C(q)\left[1+\Sigma\right]\inf_{m\in\mathcal{M}}\left\{\inf_{t\in\overline{S}_m}\|\sqrt{s_1}-t\|_2 + \sqrt{\frac{\overline{D}_m \vee \Delta_m}{n}}\right\}^q.$$

Writtten in this form, our result appears as a complete analogue of Theorem 6 of Birgé [9] about density estimation, the normalized loss function $(H/\sqrt{n})^q$ playing the role of the Hellinger loss $h^q$ for densities. We also explained in Birgé [9], Section 8.3.3, that there is a complete parallel between density estimation and estimation in the white noise model. We can therefore extend this parallel to the estimation of the intensity of a Poisson process. This parallel has also been explained and applied to various examples in Baraud and Birgé [4], Section 4.2. As an additional consequence, all the families of models that we have introduced in Sections 3.3, 4.2, 4.3 and 4.4 could be used as well for adaptive estimation of densities or in the white noise model and added to the examples given in Birgé [9].

To recover the familiar rates of convergence that we get when estimating densities which belong to some given function class $\mathcal{S}$, we merely have to assume that $s_1$ (rather than $s$) belongs to the class $\mathcal{S}$ and use the normalized loss function. Let us, for instance, apply this approach to intensities belonging to Besov spaces, assuming that $\sqrt{s_1} \in B_{p,\infty}^\alpha([0,1]^l)$ with $\alpha > l(1/p - 1/2)_+$ and that $|\sqrt{s_1}|_{B_{p,\infty}^\alpha} \leq L$ with $L > 0$. It follows that $\sqrt{s} \in B_{p,\infty}^\alpha([0,1]^l)$ with $|\sqrt{s}|_{B_{p,\infty}^\alpha} \leq L\sqrt{n}$. For $n$ large enough, $L\sqrt{n} \geq 1$ and Theorem 4 applies, leading to $\mathbb{E}_s[H^2(s,\hat{s})] \leq C(\alpha,p,l)(L\sqrt{n})^{2l/(2\alpha+l)}$. Hence

$$\mathbb{E}_s\left[n^{-1}H^2(s,\hat{s})\right] \leq C(\alpha,p,l)L^{2l/(2\alpha+l)}n^{-2\alpha/(2\alpha+l)},$$

which is exactly the result we get for density estimation with $n$ i.i.d. observations.

The same argument can be developed for the problem we considered in Section 4.2. If we assume that $\sqrt{s_1}$, rather than $\sqrt{s}$, belongs to $\mathcal{H}(\boldsymbol{\alpha},\boldsymbol{R})$, then $\sqrt{s} \in \mathcal{H}(\boldsymbol{\alpha},\sqrt{n}\boldsymbol{R})$ and the condition $R_j \geq \eta$ of Corollary 2 becomes, after this rescaling, $\sqrt{n}R_j \geq (\sqrt{n}\overline{R})^{k/(2\overline{\alpha}+k)}$ which always holds for $n$ large enough. The corresponding normalized risk bound can then be written

$$\mathbb{E}_s\left[n^{-1}H^2(s,\hat{s})\right] \leq C(k,\boldsymbol{\alpha})\overline{R}^{2k/(2\overline{\alpha}+k)}n^{-2\overline{\alpha}/(2\overline{\alpha}+k)},$$

which corresponds to the rate of convergence for this problem in density estimation.

Another interesting case is the one considered in Section 4.4. Let us assume here that instead of putting the assumptions of Proposition 5 on $\sqrt{s}$ we put them on $\sqrt{s_1}$. This implies that $\sqrt{s}$ satisfies the same assumptions with $R$ replaced by $R\sqrt{n}$ and $A$ by $A\sqrt{n}$. Then, for $n \geq n_0(A,R,\delta)$, $\gamma \leq 2\delta^{-1}\log n \leq nR^2$ and

$$\mathbb{E}_s\left[n^{-1}H^2(s,\hat{s})\right] \leq C(q,\delta,A,R)\left(n^{-1}\log n\right)^{1-q/2}.$$

This result is comparable to the bounds obtained in Corollary 3.1 of Kerkyacharian and Picard [23] but here we do not know the relationship between $q$ and $\delta$. For the special situation of $\sqrt{s_1} \in BV([0,1)^2)$, we get $\mathbb{E}_s[n^{-1}H^2(s,\hat{s})] \leq C(q,\delta,s_1) \times (n^{-1}\log n)^{1/2}$. One could also translate all other risk bounds in the same way.

An alternative asymptotic approach, which has been considered in Reynaud-Bouret [30], is to assume that $\boldsymbol{X}$ is a Poisson process on $\mathbb{R}^k$ with intensity $s$ with respect to the Lebesgue measure on $\mathbb{R}^k$, but which is only observed on $[0,T]^k$. We therefore estimate $s\mathbb{1}_{[0,T]^k}$, letting $T$ go to infinity to get an asymptotic result. We only assume that $\int_{[0,T]^k} s(x)\,dx$ is finite for all $T > 0$, not necessarily that $\int_{\mathbb{R}^k} s(x)\,dx < +\infty$. For simplicity, let us consider the case of intensities $s$ on $\mathbb{R}^+$ with $\sqrt{s}$ belonging to the Hölder class $\mathcal{H}(\alpha, R)$. For $t$ an intensity on $\mathbb{R}^+$, we set for $0 \leq x \leq 1$, $t_T(x) = Tt(Tx)$ so that $t_T$ is an intensity on $[0,1]$ and $H(t_T, u_T) = H(t\mathbb{1}_{[0,T]}, u\mathbb{1}_{[0,T]})$. Since $\sqrt{s_T} \in \mathcal{H}(\alpha, RT^{\alpha+1/2})$ it follows from Corollary 2 that there is a T-estimator $\hat{s}_T(\boldsymbol{X})$ of $s_T$ satisfying

$$\mathbb{E}_s\left[H^2(s_T, \hat{s}_T)\right] \leq C(\alpha)\left(RT^{\alpha+1/2}\right)^{2/(2\alpha+1)} = C(\alpha)TR^{2/(2\alpha+1)}.$$

Finally setting $\hat{s}(y) = T^{-1}\hat{s}_T(T^{-1}y)$ for $y \in [0,T]$, we get an estimator $\hat{s}(\boldsymbol{X})$ of $s\mathbb{1}_{[0,T]}$ depending on $T$ with the property that

$$\mathbb{E}_s\left[H^2(s\mathbb{1}_{[0,T]}, \hat{s})\right] \leq C(\alpha)TR^{2/(2\alpha+1)} \quad \text{for all } T > 0.$$

### 4.7. An illustration with Poisson regression

As we mentioned in the introduction, a particular case occurs when $\mathcal{X}$ is a finite set that we shall assume here, for simplicity, to be $\{1; \ldots; 2^n\}$. In this situation, observing $\boldsymbol{X}$ amounts to observing $N = 2^n$ independent Poisson variables with respective parameters $s_i = s(i)$ where $s$ denotes the intensity with respect to the counting measure. If we introduce a family of linear models $\overline{S}_m$ in $\mathbb{R}^N$ to approximate $\sqrt{s} \in \mathbb{R}^N$ with respect to the Euclidean distance, we simply apply Theorem 1 to get the resulting risk bounds. In this situation, the Hellinger distance between two intensities is merely the Euclidean distance between their square roots, up to a factor $1/\sqrt{2}$.

As an example, we shall consider linear models spanned by piecewise constant functions on $\mathcal{X}$ as described in Section 1.4, i.e. $\overline{S}_m = \{\sum_{j=1}^D a_j \mathbb{1}_{I_j}\}$ when $m = \{I_1, \ldots, I_D\}$ is a partition of $\mathcal{X}$ into $D = |m|$ nonvoid intervals. In order to define suitable weights $\Delta_m$, we shall distinguish between two types of partitions. First we consider the family $\mathcal{M}_{BT}$ of dyadic partitions derived from binary trees and described in Section 4.3. We already know that the choice $\Delta_m = 2|m|$ is suitable for those partitions and (4.4) applies. Note that these include the regular partitions, i.e. those for which all intervals $I_j$ have the same size $N/|m|$ and $|m| = 2^k$ for $0 \leq k \leq n$. For all other partitions, we simply set $\Delta_m = \log\binom{N}{|m|} + 2\log(|m|)$ so that (1.13) holds with $\Sigma < 3$ since the number of possible partitions of $\mathcal{X}$ into $|m|$ intervals is $\binom{N-2}{|m|-1}$. We omit the details. Denoting by $\|\cdot\|_2$ the Euclidean norm in $\mathbb{R}^N$, we derive from Theorem 1 the following risk bound for T-estimators:

$$\mathbb{E}_s\left[\left\|\sqrt{s} - \sqrt{\hat{s}}\right\|_2^2\right]$$
$$\leq C\left[\inf_{m \in \mathcal{M}_{BT}}\left\{\inf_{t \in \overline{S}_m}\|\sqrt{s} - t\|_2^2 + |m|\right\}\right.$$
$$\left.\bigwedge \inf_{m \in \mathcal{M} \setminus \mathcal{M}_{BT}}\left\{\inf_{t \in \overline{S}_m}\|\sqrt{s} - t\|_2^2 + \log(|m|) + \log\binom{N}{|m|}\right\}\right].$$

The performance of the estimator then depends on the approximation properties of the linear spaces $\overline{S}_m$ with respect to $\sqrt{s}$. For instance, if $\sqrt{s}$ varies regularly, i.e. $|\sqrt{s_i} - \sqrt{s_{i-1}}| \leq R$ for all $i$, one uses a regular partition which belongs to $\mathcal{M}_{BT}$ to approximate $\sqrt{s}$. If $\sqrt{s}$ has bounded $\alpha$-variation, as defined in Section 4.3, one uses dyadic partitions as explained in this section. If $\sqrt{s}$ is piecewise constant with $k$ jumps, it belongs to some $\overline{S}_m$ and we get a risk bound of order $\log(k+1) + \log\binom{N}{k+1}$.

## 5. Aggregation of estimators

In this section we assume that we have at our disposal a family $\{\hat{s}_m, m \in \mathcal{M}'\}$ of intensity estimators, (T-estimators or others) and that we want to select one of them or combine them in some way in order to get an improved estimator. We already explained in Section 2.3 how to use the procedure of thinning to derive from a Poisson process $\boldsymbol{X}$ with mean measure $\mu$ two independent Poisson processes with mean measure $\mu/2$. Since estimating $\mu/2$ is equivalent to estimating $\mu$, we shall assume in this section that we have at our disposal two independent processes $\boldsymbol{X_1}$ and $\boldsymbol{X_2}$ with the same unknown mean measure $\mu_s$ with intensity $s$ to be estimated. We assume that the initial estimators $\hat{s}_m(\boldsymbol{X_1})$ are all based on the first process and therefore independent of $\boldsymbol{X_2}$. Proceeding conditionally on the first process, we use the second one to mix the estimators.

We shall consider here two different ways of aggregating estimators. The first one is suitable when we want to choose one estimator in a large (possibly infinite) family of estimators and possibly attach to them different prior weights. The second method tries to find the best linear combination from a finite family of estimators of $\sqrt{s}$.

### *5.1. Estimator selection*

Here we start from a finite or countable family $\{\hat{s}_m, m \in \mathcal{M}\}$ of intensity estimators and a family of weights $\Delta_m \geq 1/10$ satisfying (1.13). Our purpose is to use the process $\boldsymbol{X_2}$ to find a close to best estimator among the family $\{\hat{s}_m(\boldsymbol{X_1}), m \in \mathcal{M}\}$.

#### *5.1.1. A general result*

Considering each estimator $\hat{s}_m(\boldsymbol{X_1})$ as a model $S_m = \{\hat{s}_m(\boldsymbol{X_1})\}$ with one single point, we set $\eta_m^2 = 84\Delta_m$. Then $S_m$ is a T-model with parameters $\eta_m, 1/2$ and $B' = e^{-2}$, (3.2) and (3.3) hold and Theorem 3 applies. Since each model is reduced to one point, one can find a selection procedure $\hat{m}(\boldsymbol{X_2})$ such that the estimator $\tilde{s}(\boldsymbol{X_1}, \boldsymbol{X_2}) = \hat{s}_{\hat{m}(\boldsymbol{X_2})}(\boldsymbol{X_1})$ satisfies the risk bound

$$\mathbb{E}_s\left[H^2(s, \tilde{s}) \big| \boldsymbol{X_1}\right] \leq C[1+\Sigma] \inf_{m \in \mathcal{M}} \left\{H^2\left(s, \hat{s}_m(\boldsymbol{X_1})\right)^2 + \Delta_m\right\}.$$

Integrating with respect to the process $\boldsymbol{X_1}$ gives

(5.1) $$\mathbb{E}_s\left[H^2(s, \tilde{s})\right] \leq C[1+\Sigma] \inf_{m \in \mathcal{M}} \left\{\mathbb{E}_s\left[H^2\left(s, \hat{s}_m\right)\right] + \Delta_m\right\}.$$

This result completely parallels the one obtained for density estimation in Section 9.1.2 of Birgé [9].

*5.1.2. Application to histograms*

The simplest estimators for the intensity $s$ of a Poisson process $\boldsymbol{X}$ are histograms. Let $m$ be a finite partition $m = \{I_1, \ldots, I_D\}$ of $\mathcal{X}$ such that $\lambda(I_j) > 0$ for all $j$. To this partition corresponds the linear space of piecewise constant functions on the partition $m$: $\overline{S}_m = \{\sum_{j=1}^D a_j \mathbb{1}_{I_j}\}$, the projection $\bar{s}_m$ of $s$ onto $\overline{S}_m$ and the corresponding histogram estimator $\hat{s}_m$ of $s$ given respectively by $\bar{s}_m = \sum_{j=1}^D (\int_{I_j} s \, d\lambda) \times [\lambda(I_j)]^{-1} \mathbb{1}_{I_j}$ and $\hat{s}_m = \sum_{j=1}^D N_j [\lambda(I_j)]^{-1} \mathbb{1}_{I_j}$ with $N_j = \sum_{i=1}^N \mathbb{1}_{I_j}(X_i)$. It is proved in Baraud and Birgé [4], Lemma 2, that $H^2(s, \bar{s}_m) \leq 2H^2(s, \overline{S}_m)$. Moreover, one can show an analogue of the risk bound obtained for the case of density estimation in Birgé and Rozenholc [13], Theorem 1. The proof is identical, replacing $h$ by $H$, $n$ by 1 and the binomial distribution of $N$ by a Poisson distribution. This leads to the risk bound

$$\mathbb{E}_s \left[H^2(s, \hat{s}_m)\right] \leq H^2(s, \bar{s}_m) + D/2 \leq 2H^2(s, \overline{S}_m) + |m|/2.$$

If we are given an arbitrary family $\mathcal{M}$ of partitions of $\mathcal{X}$ and a corresponding family of weights $\{\Delta_m, m \in \mathcal{M}\}$ satisfying (1.13) and $\Delta_m \geq |m|/2$, we may apply the previous aggregation method which will result in an estimator $\tilde{s}(\boldsymbol{X_1}, \boldsymbol{X_2}) = \hat{s}_{\hat{m}(\boldsymbol{X_2})}(\boldsymbol{X_1})$ where $\hat{m}(\boldsymbol{X_2})$ is a data-selected partition. Finally,

(5.2) $$\mathbb{E}_s \left[H^2(s, \tilde{s})\right] \leq C[1 + \Sigma] \inf_{m \in \mathcal{M}} \left\{H^2(s, \overline{S}_m) + \Delta_m\right\}.$$

Various choices of partitions and weights have been described in Baraud and Birgé [4] together with their approximation properties with respect to different classes of functions. Numerous illustrations of applications of (5.2) can therefore be found there.

## 5.2. Linear aggregation

Here we start with a finite family $\{\hat{s}_i(\boldsymbol{X_1}), 1 \leq i \leq n\}$ of intensity estimators. We choose for $\mathcal{M}$ the set of all nonvoid subsets of $\{1, \ldots, n\}$ and to each such subset $m$, we associate the $|m|$-dimensional linear subspace $\overline{S}_m$ of $\mathbb{L}_2(\lambda)$ given by

(5.3) $$\overline{S}_m = \left\{\sum_{j \in m} \lambda_j \sqrt{\hat{s}_j(\boldsymbol{X_1})} \quad \text{with } \lambda_j \in \mathbb{R} \text{ for } j \in m\right\}.$$

We then set $\Delta_m = \log \binom{n}{|m|} + 2\log(|m|)$ so that (1.13) holds with $\Sigma = \sum_{i=1}^n i^{-2}$. We may therefore apply Theorem 1 to the process $\boldsymbol{X_2}$ and this family of models conditionally to $\boldsymbol{X_1}$, which results in the bound

$$\mathbb{E}_s \left[H^2(s, \hat{s}) \big| \boldsymbol{X_1}\right] \leq C[1 + \Sigma]$$
$$\times \inf_{m \in \mathcal{M}} \left\{\inf_{t \in \overline{S}_m} \left\|\sqrt{s} - t(\boldsymbol{X_1})\right\|_2^2 + \log \binom{n}{|m|} + \log(|m|)\right\}.$$

Note that the restriction of this bound to subsets $m$ such that $|m| = 1$ corresponds to a variant of estimator selection and leads, after integration, to

$$\mathbb{E}_s \left[H^2(s, \hat{s})\right] \leq C[1 + \Sigma] \inf_{1 \leq i \leq n} \left\{\inf_{\lambda > 0} \mathbb{E}_s \left[\left\|\sqrt{s} - \lambda \sqrt{\hat{s}_i(\boldsymbol{X_1})}\right\|_2^2\right] + \log n\right\}.$$

This can be viewed as an improved version of (5.1) when we choose equal weights.

## 6. Testing balls in $(\mathcal{Q}_+(\mathcal{X}), H)$

### *6.1. The construction of robust tests*

In order to use Theorem 3, we have to find tests $\psi_{t,u}$ satisfying the conclusions of Proposition 1. These tests are provided by a straightforward corollary of the following theorem.

**Theorem 5.** *Given two elements $\pi_c$ and $\nu_c$ of $\mathcal{Q}_+(\mathcal{X})$ with respective densities $d\pi_c$ and $d\nu_c$ with respect to some dominating measure $\lambda \in \mathcal{Q}_+(\mathcal{X})$ and a number $\xi \in (0, 1/2)$, let us define $\pi_m$ and $\nu_m$ in $\mathcal{Q}_+(\mathcal{X})$ by their densities $d\pi_m$ and $d\nu_m$ with respect to $\lambda$ in the following way:*

$$\sqrt{d\pi_m} = \xi\sqrt{d\nu_c} + (1-\xi)\sqrt{d\pi_c} \quad \text{and} \quad \sqrt{d\nu_m} = \xi\sqrt{d\pi_c} + (1-\xi)\sqrt{d\nu_c}.$$

*Then for all $x \in \mathbb{R}$, $\mu \in \mathcal{Q}_+(\mathcal{X})$ and $\boldsymbol{X}$ a Poisson process with mean measure $\mu$,*

$$\mathbb{P}_\mu \left[ \log\left(\frac{dQ_{\pi_m}}{dQ_{\nu_m}}(\boldsymbol{X})\right) \geq 2x \right] \leq \exp\left[-x + (1-2\xi)\left(\frac{2}{\xi}H^2(\mu, \nu_c) - H^2(\pi_c, \nu_c)\right)\right]$$

*and*

$$\mathbb{P}_\mu \left[ \log\left(\frac{dQ_{\pi_m}}{dQ_{\nu_m}}(\boldsymbol{X})\right) \leq 2x \right] \leq \exp\left[x + (1-2\xi)\left(\frac{2}{\xi}H^2(\mu, \pi_c) - H^2(\pi_c, \nu_c)\right)\right].$$

**Corollary 6.** *Let $\pi_c$ and $\nu_c$ be two elements of $\mathcal{Q}_+(\mathcal{X})$, $0 < \xi < 1/2$ and*

$$T(\boldsymbol{X}) = \log\left((dQ_{\pi_m}/dQ_{\nu_m})(\boldsymbol{X})\right) - 2x,$$

*with $\pi_m$ and $\nu_m$ given by Theorem 5. Define a test function $\psi$ with values in $\{\pi_c, \nu_c\}$ by $\psi(\boldsymbol{X}) = \pi_c$ when $T(\boldsymbol{X}) > 0$, $\psi(\boldsymbol{X}) = \nu_c$ when $T(\boldsymbol{X}) < 0$ ($\psi(\boldsymbol{X})$ being arbitrary if $T(\boldsymbol{X}) = 0$). If $\boldsymbol{X}$ is a Poisson process with mean measure $\mu$, then*

$$\mathbb{P}_\mu[\psi(\boldsymbol{X}) = \pi_c] \leq \exp\left[-x - (1-2\xi)^2 H^2(\pi_c, \nu_c)\right] \quad \text{if } H(\mu, \nu_c) \leq \xi H(\pi_c, \nu_c)$$

*and*

$$\mathbb{P}_\mu[\psi(\boldsymbol{X}) = \nu_c] \leq \exp\left[x - (1-2\xi)^2 H^2(\pi_c, \nu_c)\right] \quad \text{if } H(\mu, \pi_c) \leq \xi H(\pi_c, \nu_c).$$

To derive Proposition 1 we simply set $\pi_c = \mu_t$, $\nu_c = \mu_u$, $\xi = 1/4$, $x = [\eta^2(t) - \eta^2(u)]/4$ and define $\psi_{t,u} = \psi$ in Corollary 6. As to (3.5), it follows from the second bound of Theorem 5.

### *6.2. Proof of Theorem 5*

It is based on the following technical lemmas.

**Lemma 8.** *Let $f$, $g$, $f' \in \mathbb{L}_2^+(\lambda)$ and $\|g/f\|_\infty \leq K$. Denoting by $\langle \cdot, \cdot \rangle$ and $\|\cdot\|_2$ the scalar product and norm in $\mathbb{L}_2(\lambda)$, we get*

(6.1) $$\int g f^{-1} f'^2 \, d\lambda \leq K\|f - f'\|_2^2 + 2\langle g, f' \rangle - \langle g, f \rangle.$$

*Proof.* Denoting by $Q$ the left-hand side of (6.1) we write

$$Q = \int gf^{-1}(f'-f)^2\, d\lambda + 2\int gf'\, d\lambda - \int gf\, d\lambda,$$

hence the result. □

**Lemma 9.** *Let $\mu, \pi$ and $\nu$ be three mean measures with $\pi \ll \nu$ and $\|d\pi/d\nu\|_\infty \leq K^2$ and let $\boldsymbol{X}$ be a Poisson process with mean measure $\mu$. Then*

$$\mathbb{E}_\mu\left[\sqrt{\frac{dQ_\pi}{dQ_\nu}}(\boldsymbol{X})\right] \leq \exp\left[2KH^2(\mu,\nu) - 2H^2(\pi,\mu) + H^2(\pi,\nu)\right].$$

*Proof.* By (1.3) and (1.2),

$$\mathbb{E}_\mu\left[\sqrt{\frac{dQ_\pi}{dQ_\nu}}(\boldsymbol{X})\right] = \exp\left[\frac{\nu(\mathcal{X}) - \pi(\mathcal{X})}{2}\right]\mathbb{E}_\mu\left[\prod_{i=1}^{N}\sqrt{\frac{d\pi}{d\nu}}(X_i)\right]$$

$$= \exp\left[\frac{\nu(\mathcal{X}) - \pi(\mathcal{X})}{2} + \int_\mathcal{X}\left(\sqrt{\frac{d\pi}{d\nu}}(x) - 1\right)d\mu(x)\right]$$

$$= \exp\left[\frac{\nu(\mathcal{X}) - \pi(\mathcal{X})}{2} - \mu(\mathcal{X}) + \int_\mathcal{X}\sqrt{\frac{d\pi}{d\nu}}(x)\, d\mu(x)\right].$$

Using Lemma 8 and (1.7), we derive that

$$\int_\mathcal{X}\sqrt{\frac{d\pi}{d\nu}}(x)\, d\mu(x) \leq 2KH^2(\mu,\nu) + 2\int\sqrt{d\pi d\mu} - \int\sqrt{d\pi d\nu}$$

$$= 2KH^2(\mu,\nu) - 2H^2(\pi,\mu) + \pi(\mathcal{X}) + \mu(\mathcal{X})$$

$$+ H^2(\pi,\nu) - (1/2)[\pi(\mathcal{X}) + \nu(\mathcal{X})].$$

The conclusion follows. □

To prove Theorem 5, we may assume (changing $\lambda$ if necessary) that $\mu \ll \lambda$ and set $v = \sqrt{d\mu/d\lambda}$. We also set $t_c = \sqrt{d\pi_c/d\lambda}$, $u_c = \sqrt{d\nu_c/d\lambda}$, $t_m = \xi u_c + (1-\xi)t_c$ and $u_m = \xi t_c + (1-\xi)u_c$. Then $\pi_m = t_m^2 \cdot \lambda$ and $\nu_m = u_m^2 \cdot \lambda$. Note that $t_c, u_c, t_m, u_m$ and $v$ belong to $\mathbb{L}_2^+(\lambda)$ and that for two elements $w, z$ in $\mathbb{L}_2^+(\lambda)$, $\|w - z\|_2^2 = 2H^2(w^2 \cdot \lambda, z^2 \cdot \lambda)$. Since $\|t_m/u_m\|_\infty \leq (1-\xi)/\xi$, we may apply Lemma 9 with $K = (1-\xi)/\xi$ to derive that

$$L = \log\left(\mathbb{E}_\mu\left[\sqrt{\frac{dQ_{\pi_m}}{dQ_{\nu_m}}}(\boldsymbol{X})\right]\right) \leq \frac{1-\xi}{\xi}\|v - u_m\|_2^2 - \|v - t_m\|_2^2 + \frac{\|t_m - u_m\|_2^2}{2}.$$

Using the fact that

$$v - u_m = v - u_c + \xi(u_c - t_c), \qquad v - t_m = v - u_c + (1-\xi)(u_c - t_c),$$

$$t_m - u_m = (1 - 2\xi)(t_c - u_c)$$

and expending the squared norms, we get, since the scalar products cancel,

$$L \leq \frac{1-2\xi}{\xi}\|v - u_c\|_2^2 + \left[\xi(1-\xi) - (1-\xi)^2 + \frac{(1-2\xi)^2}{2}\right]\|t_c - u_c\|_2^2,$$

which shows that

$$L \leq (1 - 2\xi) \left[ 2\xi^{-1} H^2(\mu, \nu_c) - H^2(\pi_c, \nu_c) \right].$$

The exponential inequality then implies that

$$\mathbb{P}_\mu \left[ \log \left( \frac{dQ_{\pi_m}}{dQ_{\nu_m}}(\boldsymbol{X}) \right) \geq 2x \right] \leq e^{-x} \mathbb{E}_\mu \left[ \sqrt{\frac{dQ_{\pi_m}}{dQ_{\nu_m}}(\boldsymbol{X})} \right] = \exp[-x + L],$$

which proves the first error bound. The second one can be proved in the same way.

## Acknowledgments

Many thanks to Philippe Bougerol for some exchanges about Poisson processes and to Albert Cohen and Ron DeVore for several illuminating discussions on approximation theory, the subtleties of the space $BV(\mathbb{R}^2)$ and adaptive approximation methods. I also would like to thank the participants of the Workshop *Asymptotics: particles, processes and inverse problems* that was held in July 2006 in Leiden for their many questions which led to various improvements of the paper.